\documentclass[12pt]{article}
\usepackage[all]{xy}
\usepackage{amssymb,amsmath,amsthm,amscd,amsfonts,MnSymbol}
\usepackage{hyperref}
\newtheorem{df}{Definition}[section]
\newtheorem{thm}[df]{Theorem}

\newtheorem{pro}[df]{Proposition}
\newtheorem{lem}[df]{Lemma}
\newtheorem{cor}[df]{Corollary}
\newtheorem{rk}[df]{Remark}
\newtheorem{ex}[df]{Example}
\newcommand{\arrow}{{\longrightarrow}}

\newcommand{\E}{{\cal E}}
\newcommand{\B}{{\cal B}}

\def\im{{\rm im}}

\def\id{{\rm id}}

\date{}
\title{\bf Weak  Topologies on Toposes}
\author{{\bf Zeinab Khanjanzadeh}\and {\bf Ali Madanshekaf}\\
    Department of Mathematics\\
Faculty of Mathematics, Statistics and Computer Science\\
Semnan University\\ Semnan\\ Iran\\ emails:
z.khanjanzadeh@gmail.com\\ \qquad\qquad amadanshekaf@semnan.ac.ir}
\date{}
\begin{document}
\maketitle
%--------------------------------------------------------------------------------------
\begin{abstract}
This paper deals with the  notion of weak Lawvere-Tierney topology on a topos.
 Our  motivation to study such a notion is based on the observation that the composition of two Lawvere-Tierney topologies on a topos is no longer idempotent,
  when seen as a closure operator. For a given  topos  $\mathcal{E}$, in this paper we investigate some properties of this notion. Among other things,
  it is shown that the set  of all weak Lawvere-Tierney topologies  on $\mathcal{E}$
constitute a complete residuated lattice provided that $\mathcal{E}$ is (co)complete. Furthermore, when the weak Lawvere-Tierney topology on $\mathcal{E}$
preserves binary meets  we give an explicit description  of the (restricted) associated sheaf functor on $\mathcal{E}$.
\end{abstract}
%-----------------------------------------------------------------------------------
AMS {\it subject classification}: 18B25; 18F20; 06F05; 20M30. \\
{\it key words}: Lawvere-Tierney topology; residuated lattice; sheaf; separated object.
%%%%%%%%%%%%%%%%%%%%%%%%%%%%%%%%%%%%%%%%%%%%%%%%%%%%%%%%%%%%%%%%%%%%%%%%%%%%%%%%%%%%%%%
\section{Introduction}
%%%%%%%%%%%%%%%%%%%%%%%%%%%%%%%%%%%%%%%%%%%%%%%%%%%%%%%%%%%%%%%%%%%%%%%%%%%%%%%%%%%%%%%
It is known~\cite{jhonstone} that the set of all Lawvere-Tierney
(LT-) topologies on an elementary topos $\mathcal{E}$ is not
closed under the formation of composition.
 In fact, the composition of two LT-topologies on $\mathcal{E}$ preserves the top element of the subobject classifier $\Omega$  and  preserves  binary
  meets too, but  fails to be idempotent. From this point of view, it is natural to study LT-topologies  without idempotency.
  In point-set  topology,  closure operators without idempotency are worthwhile~\cite{Dikranjan}.  An illuminating point of view on this notion is
   provided by the concept of a $\check{{\rm C}}$ech
closure operator (or a preclosure operator in the sense of~\cite{Claes}),
originally introduced by $\check{{\rm C}}$ech in~\cite{Cech} which
is a closure operator  without idempotency. In this approach, a
pretopological space will be a set equipped with a  $\check{{\rm
        C}}$ech closure operator~\cite{Dikranjan}.

In topos theory, an analogous notion is a weak LT-topology  (or a weak
topology, for short) on an arbitrary topos. Considering
LT-topologies  in the  framework
of~\cite{Barr}, a weak topology is exactly an
LT-topology  without idempotency. The term `weak
Lawvere-Tierney topology' was coined by Hosseini and Mousavi
in~\cite{hosseini}.  On the other hand, modal closure operators on a category and its types are of interest to some mathematicians, e.g.
see~\cite{Castellini} and \cite{Dikranjan}.  The correspondence between weak topologies and modal closure operators in a topos is given here.

We remark that ``universal closure operator" in literature has
been defined by two  ways which one can prove that these two definitions are the same. Here, we illustrate  them:\\
{\bf In the sense of Borceux in~\cite{handbook1}} : An operator on
the subobjects of each object $E$ of $\mathcal{E}$
$$A\mapsto \overline{A}:{\rm Sub}_{\mathcal{E}}(E)\rightarrow {\rm Sub}_{\mathcal{E}}(E);$$
is a {\it universal closure operator} if and only if it satisfy
the following properties, where $A, B\in {\rm
Sub}_{\mathcal{E}}(E)$ and $f : F\rightarrow E$  is a morphism of
 $\E$:\\
  (B1) $A\subseteq \overline{A}$;\\
  (B2) $A\subseteq B\Rightarrow \overline{A}\subseteq \overline{B}$; \\
  (B3) $\overline{\overline{A}} = \overline{A}$;\\
  (B4) $\overline{f^{-1}(A)} = f^{-1}(\overline{A}),$
 where $f^{-1}$ is the inverse image map.\\
 This types of operators are called {\it topology} by {\bf Barr
 in~\cite{Barr}}.\\
{\bf In the sense of Maclane in~\cite{maclane}} : An operator on
the subobjects of each object $E$ of $\mathcal{E}$
$$A\mapsto \overline{A}:{\rm Sub}_{\mathcal{E}}(E)\rightarrow {\rm Sub}_{\mathcal{E}}(E);$$
is a {\it universal closure operator} if and only if it satisfy
the following properties, where $A, B\in {\rm
Sub}_{\mathcal{E}}(E)$ and $f : F\rightarrow E$  is a morphism of
 $\E$:\\
  (M1) $A\subseteq \overline{A}$;\\
  (M2) $\overline{A\cap B} =\overline{A}\cap \overline{B}$; \\
  (M3) $\overline{\overline{A}} = \overline{A}$;\\
  (M4) $\overline{f^{-1}(A)} = f^{-1}(\overline{A}).$

 The difference between these two definitions is in the conditions (B2) and (M2).
  It is easy to see that (B2) and (B3) implies (M2). We know\cite{handbook1,maclane} that LT-topologies and universal closure operators on an arbitrary topos $\E$ are in one to one correspondence.

We start in Section 2 to study weak topologies on a topos $\mathcal{E}$ and
then some properties are investigated. Additionally, we introduce a class of weak topologies, we call it {\it weak ideal topology}, on the topos
 $M$-{\bf Sets}  by means of left ideals of the monoid $M$. Here $M$-{\bf Sets} is the topos consisting of all representations $X\times M\dot\arrow X$
 of a fixed monoid $M$ on a variable set $X$. A morphism of $M$-{\bf Sets} is a function which respects the action. Also, we show that for a productive
  weak topology $j$ on $\E$, the full subcategory  of all $j$-sheaves of $\E$, namely ${\bf Sh}_j\mathcal{E}$, is a topos. In section 3, we introduce
   and study a class of weak topologies on $\E$ induced by natural transformations of the identity functor on $\E$, i.e. $\id_\E$, and we show that
    in the special case of  the topos $M$-{\bf Sets} they  correspond to weak ideal topologies with respect to certain left ideals of $M$.
In section 4, we show that the weak topologies on a (co)complete topos constitute a complete
residuated lattice. It is well known~\cite{jhonstone} that there is no
simple formula for the join of two topologies on a given
topos $\mathcal{E}$. We shall use weak topologies as a tool for
calculating joins of topologies on a (co)complete topos. In section 5, for a productive weak topology  $j$ on $\mathcal{E},$ we  establish  a left  adjoint to the inclusion functor from the category
${\bf Sep}_j\mathcal{E}$ of all separated objects of $\mathcal{E}$ to the full subcategory $\mathcal{C}_j$ of $\E$ consisting of all
objects $E$ of $\E$ for which the closure of diagonal subobject $\bigtriangleup_E$ of $E\times E$ is closed. Moreover we will show that the former
 category ${\bf Sep}_j\mathcal{E}$ is in fact a quasitopos whenever the topos $\mathcal{E}$ is  complete, cocomplete and co-well-powered. Finally,
  section 6  is devoted to  find  the associated sheaf  to any separated object of $\mathcal{E}$ with respect to a productive weak topology   on $\mathcal{E}$. Afterwards,  we constitute a restricted  associated sheaf functor to a productive
    weak topology $j$ on       $\mathcal{E}$. Besides, throughout these sections, we provide some examples, by means of the  weak ideal topology   on the topos $M$-{\bf Sets},   to show that some constructions given by LT-topologies  on an arbitrary topos need not to  transfer exactly to  weak topologies on a topos.
%%%%%%%%%%%%%%%%%%%%%%%%%%%%%%%%%%%%%%%%%%%%%%%%%%%%%%%
\section{Weak topology}
%%%%%%%%%%%%%%%%%%%%%%%%%%%%%%%%%%%%%%%%%%%%%%%%%%%%%%%
In this section, we present the basic properties of weak
LT-topologies  on toposes. To begin, we recall the following
definition from~\cite{hosseini}.
\begin{df}\label{weak lawv. tie. topology}
{\rm A {\it weak  LT-topology} on a topos
$\mathcal{E}$ is a morphism $j : \Omega\rightarrow \Omega$ such that:\\
(i) $j\circ {\rm true} = {\rm true}$;\\
(ii) $j\circ \wedge \leq  \wedge\circ (j\times j)$,\\
in which $\wedge : \Omega\times \Omega\arrow \Omega $ is the conjunction map on $\Omega $ and $\leq$ stands for the internal order on $\Omega$ that comes from the equalizer of $\wedge$ and the first projection $pr_1$ on $\Omega$. }
\end{df}
\begin{rk} {\rm        1) It is obvious that a morphism $j
: \Omega\rightarrow \Omega$ in a topos $\mathcal{E}$ is order
preserving if and only if $j\circ \wedge \leq \wedge\circ (j\times
j)$. Moreover, that it preserves the top element means that
$j\circ {\rm true} = {\rm true}$. Thus $j : \Omega\rightarrow \Omega$ is a weak topology on a topos $\E$ if and only if $j$ is order preserving and it preserves the top element of  $\Omega.$ \\
2) From logical point of view, it is recognized that a weak LT-topology in a topos is essentially the same thing as a {\it
uniform weakening operator}, defined in \cite{Britz}, in its
internal language. Stated in the context of local set theories, a
uniform weakening operator in a local set theory $S$ is a formula
$\bullet$ with exactly one free variable of type $\Omega$
satisfying the following conditions: for any formulas $\alpha,
\beta$, (i) $\alpha\models \bullet \alpha$, (ii) if $\alpha\models
\beta$ then $\bullet\alpha\models \bullet\beta$.
We should remark that this entails that
\[\bullet(\alpha\wedge \beta)\wedge \beta \leftmodels\models\bullet (\alpha)\wedge \beta.\]}
\end{rk}
Henceforth, by a {\it weak topology} on a topos $\mathcal{E}$ we shall mean a weak LT-topology  on $\mathcal{E}$.

Recall~\cite{hosseini} that an operator on the subobjects of each object $E$ of $\mathcal{E}$
$$A\mapsto \overline{A}:{\rm Sub}_{\mathcal{E}}(E)\rightarrow {\rm Sub}_{\mathcal{E}}(E);$$
is a {\it modal closure operator} if and only if it has, for all
$A, B\in {\rm Sub}_{\mathcal{E}}(E)$, the properties $A\subseteq
\overline{A}$, $A\subseteq B$
 yields that $\overline{A}\subseteq \overline{B}$, and for each  morphism
 $f : F\rightarrow E$ in $\E$, we have $\overline{f^{-1}(A)} = f^{-1}(\overline{A}),$
 where $f^{-1}$ is the inverse image map.
 %An examination  of the  proof of Proposition V.1.1 of ~\cite{maclane}  shows that

 Any weak topology      $j$ on $\mathcal{E}$,  determines a modal closure
 operator $A\mapsto \overline{A}$ on
 the subobjects $A\rightarrowtail E$ of each object $E,$ in such a way that given $A \stackrel{\iota}{\rightarrowtail} E$, its
  $j$-closure $A$ is that subobject
 of $E$ with characteristic function $j {\rm char}(\iota)$, shown as in the diagram below
 \begin{equation}\label{modal clos}
\SelectTips{cm}{}\xymatrix{\overline{A}\ar@{>->}[dd]_{\overline{\iota}}\ar[rrr]&&&1
\ar[dd]^{{\rm true}}\\
&A\ar[r] \ar@{>->}[d]_{\iota} & 1\ar[d]^{{\rm true}}&\\
E\ar@{=}[r]&E \ar[r]_{{\rm char}(\iota)} & \Omega\ar[r]_{j}&\Omega .}
\end{equation}
For simplicity, for any subobject $\iota : A\rightarrowtail E$ we will use ${\rm char}(\iota)$ and ${\rm char}(A)$ interchangeably.

Conversely, any modal closure operator on a topos $\mathcal{E}$ always gives a unique weak topology      $j$ as
indicated in the following pullback diagram:
\begin{eqnarray}\label{uniq. weak topo.}
\SelectTips{cm}{}\xymatrix{  \overline{1}  \ar@{>->}[d] _{\overline{{\rm true}}}
\ar[r]& 1 \ar@{>->}[d]^{{\rm true}} \\ \Omega \ar[r]^{j} & \Omega}
\end{eqnarray}
i.e., $j = {\rm char} (\overline{\rm true})$.

Thus, we record the following result for further applications.
\begin{pro}\label{modal clo. ope. and weak topo.}
    On a topos $\mathcal{E}$, weak  topologies $j$ are in one to one correspondence with modal closure operators $\overline{(\cdot )}$ defined
     on ${\rm Sub}_{\mathcal{E}}(E)$ for  all $E\in\E$.
\end{pro}
{\bf  Proof.} First of all, note that with   the  isomorphism ${\rm Hom}_{\mathcal{E}}(E, \Omega)\cong {\rm Sub}_{\mathcal{E}}(E)$ which is
 natural in $E$ we have the following commutative diagram:
\begin{eqnarray}\label{natural meet}
\SelectTips{cm}{}\xymatrix{  {\rm Sub}_\E(E) \times {\rm Sub}_\E(E)   \ar[d] ^{\cong}
    \ar[r]^{~~~~\bigcap}& {\rm Sub}_\E(E) \ar[dd]^{\cong} \\ {\rm Hom}_\E(E, \Omega)\times {\rm Hom}_\E(E, \Omega) \ar[d]^{\cong}& \\ {\rm Hom}_\E(E, \Omega\times  \Omega) \ar[r]^{~~\bigwedge_E}& {\rm Hom}_\E(E, \Omega)}
\end{eqnarray}
in which the map $\wedge_E$ is the map compose via $\wedge$ on $\Omega$. Let us now $j$ be a weak topology on $\E$, $E$  an object of $\E$ and $A$ a subobject of $E$. First, since $j\circ {\rm true} = {\rm true}$,  by diagram~(\ref{modal clos})  it is clear that  $A\subseteq \overline{A}$ as subobjects of $E$.   To prove the monotonicity of $\overline{(\cdot)}$, it suffices to show that  ${\rm char}(\overline{A\cap B}) \leq {\rm char}(\overline{A}\cap \overline{B})$ and then  we have $\overline{A\cap B} \subseteq \overline{A}\cap \overline{B}$ as subobjects of $E$, for all subobjects $A$ and $B$ of $E$.
To show that non-equality,  we have
$$\begin{array}{rcl}
{\rm char}(\overline{A\cap B})  & =& j{\rm char}(A\cap B) \\
&=&j\wedge ({\rm char}(A), {\rm char}(B)  )\\
&\leq & \wedge (j\times j) ({\rm char}(A), {\rm char}(B)  )~~~~(\textrm{by~Definition}~\ref{weak lawv. tie. topology} \textrm{(ii))}\\
&=& \wedge ({\rm char}(\overline{A}), {\rm char}(\overline{B})  ) \\
&=& {\rm char}(\overline{A}\cap \overline{B})
\end{array}$$
in which the last equality holds by the commutativity of  diagram~(\ref{natural meet}).\\
Next let  $f : F\rightarrow E$ be  a morphism in $\E$ and $A$ a
subobject  of $E$. By the Pullback Lemma it is evident that the
large  rectangle  in the diagram  below is   pullback
\begin{eqnarray*}
    \SelectTips{cm}{}\xymatrix{  f^{-1}(A)  \ar@{>->}[d]
        \ar[r]& A \ar@{>->}[d]\ar[r]&1 \ar@{>->}[d]^{{\rm true}}\\ F \ar[r]_{f} & E\ar[r]_{{\rm char}(A)}&\Omega}
\end{eqnarray*}
so we have ${\rm char}( f^{-1}(A)) =  {\rm char}(A) f$. Hence, we get
$${\rm char}(\overline { f^{-1}(A)}) = j {\rm char}( f^{-1}(A)) = j {\rm char}(A) f =  {\rm char}(\overline {A}) f = {\rm char}(f^{-1}(\overline { A})).$$
Then, $\overline { f^{-1}(A)} =  f^{-1}(\overline {A})$ as subobjects of $F$. Thus, $\overline{(\cdot)}$ is a modal closure operator.

Conversely, let $\overline{(\cdot)}$ be a modal closure operator defined  on all ${\rm Sub}_{\mathcal{E}}(E)$. Since $1=\overline {1}$ as subobjects of $\Omega$, then
${\rm char}(\bar{1}) = {\rm char}(1)$ and so $j\circ {\rm true} = {\rm true}$.
By monotonicity of $\overline{(\cdot)}$, it is straightforward to see that $\overline{A\cap B}\subseteq \overline{A}\cap \overline{B}$  for all subobjects $A$ and $B$ of $E$. Therefore, by the definition of the order on $\Omega$, ${\rm char}(\overline{A\cap B}) \leq {\rm char}(\overline{A}\cap \overline{B})$ and so,
$$\begin{array}{rcl}
j\wedge ({\rm char}(A), {\rm char}(B)  ) & =& j{\rm char}(A\cap B) \\
&=&{\rm char}(\overline{A\cap B}) \\
&\leq & {\rm char}(\overline{A}\cap \overline{B})\\
&=& \wedge (j\times j) ({\rm char}(A), {\rm char}(B))
\end{array}$$
in which the last equality holds by the commutativity of  diagram~(\ref{natural meet}). Then, by choosing $A= B= 1$ as subobjects of $\Omega$, we get  $j\circ \wedge \leq  \wedge\circ (j\times j)$. Now it is easily seen that the two constructions are mutually inverse.
$\qquad\square$

As usual, a monomorphism $k: B\rightarrowtail A$ in $\mathcal{E}$ is {\it $j$-dense} whenever
$\overline{B} = A$, and  {\it $j$-closed} if
  $\overline{B} = B$,  as subobjects of $A$.
It is convenient to see that any modal  closure operator $\overline{(\cdot )}$
  is {\it weakly hereditary}; that is the
subobject $A\subseteq \overline{A}$ is $j$-dense for each
subobject $A\rightarrowtail E$ (see also~\cite{Dikranjan}).

Moreover, an object $C$ is called  {\it $j$-sheaf} whenever for
any $j$-dense monomorphism $m : B\rightarrowtail A$, one can
uniquely  extend any morphism $h:B\rightarrow C$ in $\E$ to $A$ as
follows
$$\SelectTips{cm}{}\xymatrix{ B \ar[r]^h\ar@{>->}[d]_m & C\\  A \ar@{-->}[ur]_{g}}$$
%We say that $C$ is {\it $j$-separated} if the morphism $g$ exists,
%it is unique.
in the other words, the map ${\rm Hom}_{\mathcal{E}}(B , C)\stackrel{\overline{m}}{\longrightarrow} {\rm Hom}_{\mathcal{E}}(A , C)$
given by $g\mapsto gm$ is a bijection.

We say that $C$ is {\it $j$-separated} if in every commutative
diagram
$$
\SelectTips{cm}{} \xymatrix{B\ar@{>->}[r]^{m} &
A\ar@<1.0ex>[r]^{g} \ar@<-0.25ex>[r]_{g'} & C,}
$$
 where $m : B\rightarrowtail A$
 is a $j$-dense
monomorphism, it is already $g = g'$.

For a weak topology $j$ on $\E$, we denote the full subcategories of all $j$-sheaves and $j$-separated objects in $\E$ by ${\bf Sh}_j\mathcal{E}$ and ${\bf Sep}_j\mathcal{E}$, respectively.
%We will finish this section by the following remark:
    %Let $j$ be a weak topology on a topos $\E$.

For a weak topology   $j$ on a topos $\mathcal{E}$, the morphisms
$j, \id_{\Omega} : \Omega\rightarrow \Omega$ have an equalizer
$\Omega_j$,
\begin{eqnarray}\label{equalizer}
\SelectTips{cm}{} \xymatrix{\Omega_j\ar@{>->}[r]^{m} &
\Omega\ar@<1.0ex>[r]^{j} \ar@<-0.25ex>[r]_{\id_{\Omega}} & \Omega.}
\end{eqnarray}
 Also in a similar vein to~\cite[Lemma V.2.2]{maclane} the object $\Omega_j$ as in (\ref{equalizer})
classifies closed subobjects, in the sense that, for each object $E$ of $\mathcal{E}$, there is a bijection
\begin{equation}\label{omega j classifies closed subobjects}
{\rm Hom}_{\mathcal{E}}(E, \Omega_j)\stackrel{\sim}{\longrightarrow} {\rm ClSub}_{\mathcal{E}}(E);
\end{equation}
which is natural in $E$. Here  ${\rm ClSub}_{\mathcal{E}}(E)$ is the set  of all closed subobjects of $E$.
\begin{rk}\label{Barr idem}
\rm{That a weak topology $j$ on a topos
    $\mathcal{E}$ is idempotent if and only if it is an LT-topology on
    $\mathcal{E}$  is obtained in relating to closure
    operator corresponding to an  LT-topology  in the
    sense of Borceux or Barr which can be found in~\cite[Vol. I, p. 227]{handbook1} or~\cite{Barr}.}
\end{rk}
Let us denote the image of the weak topology $j$ by $\im (j)$ (for details,
see~\cite[p. 184]{maclane}). The following indicates that it is not always
necessary to have $\Omega_j = \im (j)$. Hence, we are unable
to construct the associated sheaf functor to a weak topology $j$ on
$\mathcal{E}$ as in~\cite{maclane}.

\begin{pro}\label{im j and idempotent}
For a weak topology $j$ on a topos $\mathcal{E}$, we factor  $j$ through its image as $\Omega\stackrel{r}{\twoheadrightarrow}
 \im (j)\stackrel{k}{\rightarrowtail}
\Omega.$ Then $j$ is an idempotent {\rm (}or equivalently,  is an
LT-topology  on  $\mathcal{E}${\rm )} if and only if
$\Omega_j = \im (j)$, as  subobjects of $\Omega$.
\end{pro}
{\bf  Proof.} {\it Necessity.} We show that $\Omega_j = \im (j)$. On the one hand, we always have $\Omega_j \subseteq \im (j)$. To observe this fact, we have
$$k(rm) = (kr)m = jm = m.$$
Hence, the composition morphism
$\Omega_j\stackrel{m}{\rightarrowtail}
\Omega\stackrel{r}{\twoheadrightarrow}\im (j)$ is the required map
which completes the triangle below.
$$\SelectTips{cm}{}\xymatrix{\Omega_j\ar@{>->}[dr]_{m} \ar@{-->}[rr]^{rm} &&\im (j)\ar@{>->}[dl]^{k}\\
&\Omega &}$$
The equality $k(rm) = m$ indicates that $rm$ is monic.
On the other hand, we establish $\im (j) \subseteq \Omega_j$. By the assumption $jj= j$,
 in the equalizer diagram (\ref{equalizer}), there exists
a unique morphism $n :\Omega\rightarrow\Omega_j$ for which $mn=
j$. This implies that $mn kr= jj= j = kr.$ Since  $r$ is an
epimorphism, $mn k= k$. Therefore, the composition morphism $\im (j)\stackrel{k}{\rightarrowtail}
 \Omega\stackrel{n}{\rightarrow}\Omega_j$ is the required map which completes the triangle below
 $$\SelectTips{cm}{}\xymatrix{\im (j)\ar@{>->}[dr]_{k} \ar@{-->}[rr]^{nk} &&\Omega_j\ar@{>->}[dl]^{m}   \\ &\Omega . &}$$

{\it Sufficiency.} Suppose that  $\Omega_j = \im (j)$. Then, the map $j$  is idempotent because:
%We show that $j$ is idempotent as follows:
$$\begin{array}{rcl}
jj & =& jkr\\
&=& jmr~~~~(k=m, \ \textrm{by\ the\ assumption})\\
&=& mr~~~~~~~~~~~~(\textrm{by}~(\ref{equalizer}))\\
&=& kr\\
&=&j
\end{array}$$
This completes the proof. $\qquad\square$

Notice that the proof of Proposition \ref{im j and idempotent}
shows that  a weak topology $j$ on $\mathcal{E}$  is idempotent if
and only if we just have $\im (j) \subseteq \Omega_j$, as  subobjects
of $\Omega$. Meanwhile, we point out that a weak topology $j$ on
$\mathcal{E}$  can
be rewritten as an LT-topology  (in the sense of
universal closure operator of Borceux or in the sense of topology of Barr) exactly  without idempotency of $j$ or explicitly
 without the condition $j\geq j^2$:
% as one has $j\circ {\rm true} = {\rm true}$ for a weak topology  $j$ on $\mathcal{E}$,
it yields from definition~\ref{weak lawv. tie. topology}(i)  that $j\leq j^2$.

%In the following example we give some examples of LT-topologies on
%toposes which the compositions of them are not an LT-topology.
At this moment, it would be good to have  some examples of the composite of two LT-topologies which is not an LT-topology.
\begin{ex}
    {\rm  Let $u: 1 \rightarrow \Omega$ be a global element of
        $\Omega$ in an arbitrary topos $\E$. By~\cite[Vol III, Example
        9.1.5.c]{handbook1}, we have two LT-topologies $u\vee -$ and
        $u\Rightarrow -$ on $\E$. In this way, the compositions of each two of three
        LT-topologies $u\vee -$, $u\Rightarrow -$ and the double negation LT-topology $\neg\neg$ on $\E$
        are not necessarily commutative and then hence are not LT-topology.\\
       }
\end{ex}
Next we provide another difference between weak topologies and  LT-topologies  on a topos.
\begin{rk}\label{not productive}
{\rm    For a weak topology      $j$ on $\mathcal{E}$, it can be easily checked that
$j\circ \wedge = {\rm char}(\overline{{\rm true}
 \times {\rm true}})$ and $\wedge\circ (j\times j)={\rm char}(\overline{{\rm true}} \times \overline{{\rm true}})$ (see~\cite{jhonstone}).
 By Definition~\ref{weak lawv. tie. topology} (ii), it follows that
two subobjects $\overline{1\times 1}$ and $\overline{1}\times \overline{1}$ of $\Omega
 \times \Omega$, are not necessarily equal. This means that the modal closure operator associated to $j$, namely $\overline{(\cdot )}$,
 is not {\it productive}; that is the closure
 does not commute with existing products in  $\mathcal{E}$ (see also,~\cite{Dikranjan}).
 It is easy to see that for a weak topology     $j$ on $\mathcal{E}$,   the modal closure operator
  associated  to $j$, i.e. $\overline{(\cdot )}$,  is productive if and only if one has $j\circ \wedge = \wedge\circ (j\times j)$ if and only if  $\overline{(\cdot )}$    commutes with binary intersections
   (see, also the proof of~\cite[Vol. III, Proposition 9.1.3]{handbook1}).}
\end{rk}
Following Remark~\ref{not productive},  a weak  topology $j$ on
$\mathcal{E}$ is called  {\it productive} if $j\circ \wedge = \wedge\circ (j\times j)$. Notice that any LT-topology  is productive. Let us now give other  examples of productive weak topologies.
\begin{ex}\label{weak ideal topology}
 {\rm  (i) The composite of any two topologies on a topos $\mathcal{E}$
is a productive weak topology (see also, \cite[Lemma A. 4.5.17]{jhonstone}).

(ii) Let  us recall that the commutative monoid of natural
endomorphisms of the identity functor on  a topos $\E$ is called
{\it the center of $\E$}.
Suppose that $\alpha$ is a natural endomorphism of the
identity functor on $\E$.  It is easy to see that
$\alpha_{\Omega}$ is a  productive weak topology on $\E$. (See also,~\cite[Lemma 1.1]{jhonstone1}).}
\end{ex}

Let us now take a monoid $M$ and a left ideal $I$ of $M$. The {\it
ideal {\rm(or} residual{\rm)} closure operator} with respect to
$I$, which can be found  in~\cite{Ebrahimi}, is given by
\begin{equation}\label{ideal closure}
\overline{A} = \{b\in E~|~\forall n\in I, ~ bn\in A\},
\end{equation}
for any sub $M$-set $A$ of an $M$-set $E$ (see also,~\cite{bir}). \\
Notice that the subobject classifier $\Omega_M$ on  $M$-{\bf Sets}
is the set of all right ideals
$K$ of $M$ endowed with the action $\cdot : \Omega_M\times M\rightarrow \Omega_M$ defined by
\begin{equation}
K\cdot m = \{n\in M \mid mn\in K\}
\end{equation}
for $K\in \Omega_M$ and $m\in M$.

Now, in the following theorem we introduce a class of (productive) weak topologies in the topos $M$-{\bf Sets}.
\begin{thm}\label{weak ideal topologies}
    Let $I$ be a left ideal of $M$. Then the  ideal closure operator $\overline{(\cdot )}$ as
    in~(\ref{ideal closure}) is a modal closure   operator on the topos  $M$-{\bf Sets}, and the weak topology  corresponding
     to the ideal  closure operator is the  action preserving  map $j^I :
    \Omega_M\rightarrow \Omega_M$ given by
    \begin{equation}\label{def. of jI}
    j^I (K) = \{m\in M\mid \forall n\in I, ~ mn\in K\}
    \end{equation}
    for any $K\in \Omega_M$.   Moreover, for a (left) ideal $I$ of $M$ the operator  $j^I$  is idempotent if and only if  {\rm (}$(IM)^2 = IM${\rm )} $I^2 = I$.
\end{thm}
{\bf  Proof.}
It is straightforward to see that $\overline{(\cdot )}$ is a modal closure   operator on the topos $M$-{\bf Sets}. Also, by the
diagram~(\ref{uniq. weak topo.}) and the identity (\ref{ideal closure})  we observe  that for any right ideal $K$ of the monoid $M$ we have
\[j^I (K) = \{m\in M \mid K\cdot m\in \overline{\{M\}}\} = \{ m\in M\mid \forall n\in I, ~ mn\in K\}.\]
Meanwhile, using~\cite[Theorem 2.10]{Ebrahimi} for any left ideal $I$ of $M$,  it is seen that $j^I$ is idempotent  if and only if $(IM)^2 = IM$ and  moreover if  $I$ is an ideal of $M$, then $j^I$ is idempotent  if and only if $I^2 = I$.   $\quad\quad\square$

By a {\it weak ideal topology} induced by the left ideal $I$ of $M$ we mean  the weak topology  $j^I$ defined in Theorem~\ref{weak ideal topologies}.
One can easily check  that $j^I$ is a
productive weak topology      on  $M$-{\bf Sets}, equivalently,
$\overline{A\cap B} = \overline{A}\cap\overline{B}$ for each $M$-set $E$ and any two sub $M$-sets $A$ and $B$ of $E$. Likewise, it is evident that a sub $M$-set $A\subseteq E$ is $j^I$-dense in  $M$-{\bf Sets} if and only if  $EI\subseteq A$.
Finally, we note that the {\it weak Grothendieck topology}  induced by the left ideal $I$ of  the monoid $M$
(as a category with just one object)
associated to $\overline{(\cdot )}$, which is defined as in~\cite{hosseini} in the general case, stands for $\mathbf{J}^I = \{K\in \Omega_M|~ I\subseteq K\}.$ For a comprehensive study of weak ideal topology on the topos $M$-{\bf Sets} we refer the reader to~\cite{Khanjanzadeh2, Khanjanzadeh}.

Now, we can give other examples of LT-topologies which the composition of them is not an LT-topology. For two  ideals $I$ and $J$ of $M$, it is clear that the product  $IJ = \{mn \mid m\in I, n\in J\}$ is an ideal of $M$.  By~(\ref{def. of jI}), it is easy to see that $j^I\circ j^J        = j^{IJ}$ and $j^J\circ j^I = j^{JI}$. On the other hand, we have $j^{IJ} = j^{JI}$  if and only if  $IJ =JI$. Thus, for two idempotent ideals  $I$ and $J$
        of $M$ with the property $IJ\not= JI$, we have $j^I\circ j^J \not
        = j^J\circ j^I$. This shows that $j^I\circ j^J$ is not an LT-topology.

In an analogous manner to LT-topologies  as in~\cite[Lemma V. 2.1]{maclane}, one can observe that
 both subcategories ${\bf Sh}_j\mathcal{E}$ and  ${\bf Sep}_j\mathcal{E}$ of $\mathcal{E}$, associated to any weak topology
  $j$  on $\mathcal{E}$, are closed under all limits which exist in $\E$, as well as under exponentiation
  with an arbitrary object of $\mathcal{E}$.
  %\marginpar{this sentence added}
  Therefore, the inclusion
  functor ${\bf Sh}_j\mathcal{E} \hookrightarrow \mathcal{E}$
  preserves finite limits and exponentials.

In the  following  lemma we show that  for a productive weak
topology $j$, $\Omega_j$ is actually  a sheaf and hence  ${\bf Sh}_j\mathcal{E}$ is a topos.
\begin{lem}\label{being j- separated of omega-j}
Let  $j$ be a productive weak    topology on $\mathcal{E}$. Then
the object $\Omega_j$ defined as in (\ref{equalizer}) is  a $j$-sheaf.
\end{lem}
{\bf  Proof.} In view of the isomorphism (\ref{omega j classifies closed subobjects}) and by the definition of a sheaf,
it is sufficient to prove that for any $j$-dense monomorphism
 $n : A\rightarrowtail E$  the inverse image map
$$n^{-1} :  {\rm ClSub}_{\mathcal{E}}(E)\longrightarrow {\rm ClSub}_{\mathcal{E}}(A)$$
 is an isomorphism. Let us first prove that $n^{-1}$ is well-defined. Suppose that
 $\iota : P\rightarrowtail E$ be a closed subobject
  of $E$. Since  $n^{-1}(P) = n^{-1}(\overline{P}) = \overline{n^{-1}(P)}$ it follows that the subobject
 $n^{-1}(P) = P\cap A$ of $A$  is $j$-closed.

Also we should remark that for two given  (not necessary closed) subobjects $P ,Q$ of $E$, if $n^{-1}(P) = n^{-1}(Q)$ then we have
$$\overline{P} =  \overline{P}\cap E= \overline{P}\cap \overline{A} =\overline{P\cap A} = \overline{Q\cap A} = \overline{Q}\cap \overline{A}
 = \overline{Q}\cap E = \overline{Q},$$ in ${\rm Sub}_{\mathcal{E}}(E)$.

Next, given a closed subobject $\iota : P\rightarrowtail A$ of $A$,
we form the composition subobject $n\iota: P\rightarrowtail E$ of
 $E$ which is the subobject $\exists_n(P)\rightarrowtail E$ (see the definition of $\exists_n(P)$ as in~\cite{maclane}). One can observe that
\begin{eqnarray}\label{n -1 p bar}
n^{-1}(\exists_n(P))  & =& P~~~~~~~~~~~~~~~~{\rm in \ {\rm Sub}_{\mathcal{E}}(A)}\nonumber\\
&=& \overline{P} ~~~~~~~~~~~~~~~~{\rm in \ {\rm Sub}_{\mathcal{E}}(A)}\nonumber\\
 & =& \overline{n^{-1}(\exists_n(P))} = n^{-1}(\overline{\exists_n(P)}),
\end{eqnarray}
in which the first and  the third equalities hold for $P\subseteq
A$ and the second one as $P$ is closed in $A$. Directly, by
(\ref{n -1 p bar}) and preceding remark  we have
$\overline{\exists_n(P)} = \overline{\overline{\exists_n(P)}}$ in
 ${\rm Sub}_{\mathcal{E}}(E),$  so $\overline{\exists_n(P)}$ is a closed subobject of $E$ such that  $ n^{-1}(\overline{\exists_n(P)}) = P$.

 On the other hand, for any  closed subobject $Q\rightarrowtail E$ of $E$, we get
\begin{eqnarray*}
\overline{\exists_nn^{-1}(Q)}  & =&  \overline{\exists_n(Q\cap A)}\\
&=& \overline{Q\cap A} ~~~~~~~~~~~~~~~~{\rm in \ {\rm Sub}_{\mathcal{E}}(E)}\nonumber\\
 & =& \overline{Q}\cap \overline{A} = Q
\end{eqnarray*}
 in which the last equality holds as $Q$ is closed in $E$ and $A$ is dense there.
 Therefore, the map $n_1 : {\rm ClSub}_{\mathcal{E}}(A)
 \longrightarrow {\rm ClSub}_{\mathcal{E}}(E)$ defined by
 $n_1(P) = \overline{\exists_n(P)}$, for any  closed subobject $P\rightarrowtail A$
  of $A$, is the inverse of $n^{-1}$. We are done. $\qquad\square$

By the above Lemma we can deduce the following Theorem.
\begin{thm}
For  a productive weak  topology $j$ on $\E$ the  category  ${\bf Sh}_j\mathcal{E}$ is a topos. In Particular, for any  left ideal $I$ of a monoid $M$, the sheaf category ${\bf Sh}_{j^I}(M\!-\bf{Sets})$ is a topos.
\end{thm}
Take a weak topology $j$ on a topos $\E$.  Let us consider in $\E$
the equalizer $w: L_j\rightarrowtail\Omega$ of the two
morphisms $j^2 , j: \Omega \rightarrow \Omega $.  For all objects
$E$ of $\E$, we shall write $ {\rm CClSub}_{\mathcal{E}}(E)$ to
denote the set of all
 subobjects $A$ of $E$ for which $\overline{A}$ is closed in  $E$. We remark that since $j$ is not idempotent the two sets $ {\rm CClSub}_{\mathcal{E}}(E)$ and $ {\rm ClSub}_{\mathcal{E}}(E)$ are not the same in general (see also Example~\ref{cou. example about sepa.}). In an analogous manner to
(\ref{omega j classifies closed subobjects})  one can observe, for all $E$ of  $\E$, the bijection
\begin{equation}\label{Lj and double closed subobjects}
{\rm Hom}_{\mathcal{E}}(E, L_j)\stackrel{\sim}{\longrightarrow} {\rm CClSub}_{\mathcal{E}}(E);
\end{equation}
which is natural in $E$. It is evident that $\Omega_j\subseteq
L_j$ as subobjects of $\Omega$. Meanwhile, the morphism ${\rm
true}:1\rightarrow\Omega$ factors through $L_j$.

The following indicates a treatment of $L_j$ which is close to
$\mathcal{M}$-injectivity property of $\Omega_j$, where $\mathcal{M}$ is the
 class of all $j$-dense monomorphisms in $\mathcal{E}$.
\begin{rk}
{\rm Let $j$ be a productive weak topology on  $\E, n : A\rightarrowtail E$
a $j$-dense monomorphism and $P\rightarrowtail A$ an
element of ${\rm CClSub}_{\mathcal{E}}(A)$.
We  denote here  the closure of $P$ as a subobject of $A$ by $C(P)$, instead of $\overline{P}$,  which  is closed in $A$, i.e. $C(C(P)) = C(P)$ by the definition of ${\rm
    CClSub}_{\mathcal{E}}(A)$. One can
 consider the inverse image map
\[ n^{-1} :  {\rm CClSub}_{\mathcal{E}}(E)\longrightarrow {\rm
CClSub}_{\mathcal{E}}(A)\]
which can be easily checked that it is
well defined. In a completely similar way as in Lemma \ref{being
j- separated of omega-j}, by replacing $C(P)$ by $P$ in that proof, we can deduce
that $\exists_n (C(P))$ belongs to ${\rm CClSub}_{\mathcal{E}}(E)$
and then, $n^{-1}(\exists_n (C(P))) = C(P)$ in
 ${\rm CClSub}_{\mathcal{E}}(A)$.}
\end{rk}
The following example shows  in  general that the object $L_j$  is not  separated.
\begin{ex}\label{cou. example about sepa.}
{\rm Take a zero semigroup $S$, that is $S\not = 0$ and $S^2 = 0$,
and the induced monoid $M= S \dot{\cup} \{1\}$ by setting $s1 = s
= 1s$ and $1\cdot 1 = 1$. Evidently, $S$ is a (two sided) ideal of $M$. Let
$j^S$ be  the weak ideal  topology with respect to the ideal $S$
of $M$ as in Theorem~\ref{weak ideal topologies}. Under these
circumstances, it is easily seen that $L_{j^S}= {\rm equal}(j^S ,
j^Sj^S) = \{\emptyset , S , M\}.$ Notice that the $M$-set $L_{j^S}$ is not
separated. To obtain  this, consider the $j^S$-dense monomorphism
$SM\subseteq M$ and two  action preserving  maps $f , g :
M\rightarrow L_{j^S}$ given by $f(1) = M$ and $g(1) = S$. From
 $M\not = S$ it follows that $f\not = g$.  But one can easily observe that the two
restriction maps
\[f|_{SM}, g|_{SM}: SM\rightarrow L_{j^S}\] are equal to the constant map in $M$. \\
We point out that the ideal $S$ of $M$ belongs to  $L_{j^S}$ since $j^Sj^S(S) = M = j^S(S)$ by the definition of $j^S$.
 Also, for the sub $M$-set $S$ of $M$,
from $\overline{\overline{S}}= M = \overline{S}$ we deduce that the sub $M$-set $S\subseteq M$ belongs to
 ${\rm CClSub}_{M-{\bf Sets}}(M)$. However, one has
$S\not = \overline{S}$. This shows that ${\rm CClSub}_{M-{\bf Sets}}(M)\supsetneqq  {\rm ClSub}_{M-{\bf Sets}}(M).$}
\end{ex}
Let $j$ be  a weak topology on a topos $\mathcal{E}$ and let $m$ and $n$ be two composable subobjects  in
 $\mathcal{E}$.  It is easy  to see that if the  composite subobject $mn$ is dense then so are $m$ and $n$.
 In contrast with Lawvere-Tierney topologies~\cite[A.4.5.11(iii)]{jhonstone}, the converse is not necessarily true.
  For example, consider a left ideal $I$ of a monoid  $M$ for which $(IM)^2\not = IM$ and the chain $I^2M\subseteq IM\subseteq M$   of  $j^I$-dense monomorphisms in $M$-{\bf Sets}. One can easily check that $I^2M\subseteq M$ is not  $j^I$-dense.
%%%%%%%%%%%%%%%%%%%%%%%%%%%%%%%%%%%%%%%%%%%%%%%%%%%%%%%%%%%%%%%%%%%%%%%%%%%%%%
\section{Natural endomorphisms and principal weak topologies}
%%%%%%%%%%%%%%%%%%%%%%%%%%%%%%%%%%%%%%%%%%%%%%%%%%%%%%%%%%%%%%%%%%%%%%%%%%%%%%
In this  section, we present  some properties of the productive
weak topology which we have already  introduced in
Example~\ref{weak ideal topology} (ii). To begin, we proceed to
extend~\cite[Definition 7.3.3]{LB} to weak topologies on  a topos $\E$.
\begin{df}\label{princi. topo.}
{\rm Let $j$ be a weak topology on  $\E$, $E$  an object and
$V\mapsto \overline{V}$ stands for the  closure operator
associated with $j$ on subobjects $V \subseteq E$. We say $j$ is
{\it principal} if, for all objects $E$, the closure operator on
${\rm Sub}_{\E}(E)$ has a left adjoint $U\mapsto U^{\circ}$,
called {\it interior}, that is,
$$U^{\circ}\subseteq V \Longleftrightarrow U\subseteq \overline{V}~ \textrm{in}~{\rm Sub}_{\E}(E).$$ }
\end{df}
In a  similar way to~\cite{LB} one gets:
\begin{lem}\label{weak. pri. and dense sub.}
A weak topology $j$ on $\E$ is  principal if and only if for all
objects $E$ of  $\E$ there exists a least dense subobject $U_E$ of $E$.
\end{lem}
We remark that, by the proof of~\cite[Lemma 7.3.5]{LB}, for any
subobject $A\rightarrowtail E$ of $\E$, indeed one has $U_A =
A^\circ$ as subobjects of $E$. Also, corresponding to any  principal
weak topology $j$, in a similar way to~\cite[p. 147]{LB},  there
exists a functor $E\mapsto E^\circ : \E\rightarrow \E$ denoted by
$U$ which we call it the {\it interior functor}.
In an analogous way to~\cite{LB}, it is immediate that:
\begin{rk}
{\rm Let $j$ be a principal weak topology on a topos $\E$. One can
easily check that:

(i) For any $E\in \E,$ $E$ is open, i.e. $E =
E^\circ$, if and only if $(\bigtriangleup_E)^\circ =
\bigtriangleup_E$ in ${\rm Sub}_{\E}(E\times E),$ in which $\bigtriangleup_E$ is the  diagonal map on $E$.

(ii) The interior functor $E\mapsto E^\circ$ preserves monos.

(iii) For any $f : X \rightarrow Y$ in $\E$ and $V \in {\rm Sub}(X)$, one has
$\exists_f (V^\circ) \cong (\exists_fV )^\circ$.

(iv) The interior functor $E\mapsto E^\circ$ is right adjoint to the
inclusion of the category of open objects of $\E,$ into
$\E$.}
\end{rk}
 For the weak topology $\alpha_{\Omega}$ on  $\E$, given in Example~\ref{weak ideal topology} (i),  a monomorphism
$m : A\rightarrowtail E$ is  $\alpha_{\Omega}$-dense  if and only if  $\im
(\alpha_E)\subseteq A,$ as subobjects of $E$. Because, $m : A\rightarrowtail E$ is  $\alpha_{\Omega}$-dense  if and only if   $\alpha_{\Omega}{\rm char}(A) = \rm{true}_E$  if and only if  ${\rm char}(A)\alpha_E = \rm{true}_E$  if and only if  there exists a unique morphism $f : E\rightarrow A$ such that $mf = \alpha_E$  if and only if  $\im
(\alpha_E)\subseteq A,$ as subobjects of $E$.

Consequently, by Lemma~\ref{weak. pri. and dense sub.}, we have the following example of principal weak topology.
\begin{pro}\label{principal of aloha}
Let $\alpha$ be a natural endomorphism of the identity functor on
$\E$. Then the weak topology $\alpha_{\Omega}$ on  $\E$ is
principal in which for all objects $E$ of  $\E$, $U_E$ is the
subobject $\im (\alpha_E)$  of $E$.
\end{pro}
Let $j$ be a weak topology on $\E$. In a similar way to \cite[Vol.
I, p. 235]{handbook1}, we will say a morphism $f : A\rightarrow B$
in $\E$ is $j$-{\it bidense} when its image is $j$-dense  and the
equalizer of its kernel pair is $j$-dense. Then a monomorphism in
$\E$ is dense if and only if it is bidense. The proof is trivial because the equalizer of a kernel pair of a mono is the
largest subobject.
%\cite[Proposition 5.8.2]{handbook1}.

Now we record the following theorem:
\begin{thm}
    Let $\alpha$ be a natural endomorphism of the identity functor on
    $\E$. Then  for all objects $E$ of $\E$,  $\alpha_E : E\rightarrow
    E$ is bidense for the weak topology $\alpha_\Omega$ on  $\E$.
\end{thm}
{\bf  Proof.} By the paragraph before Proposition \ref{principal of aloha}, it is clear that $\alpha_E$ is $\alpha_\Omega$-dense. Now we construct the kernel pair of $\alpha_E$ which is the following pullback diagram
\begin{eqnarray}\label{pulb. of alphaE}
\SelectTips{cm}{}\xymatrix{ P  \ar[d] _{v} \ar[r]^{u}&
    E \ar[d]^{\alpha_E} \\ E \ar[r]_{\alpha_E } & E.  }
\end{eqnarray}
We denote the equalizer of $u$ and $v$ by $e : B\rightarrowtail
P$. By the naturality of $\alpha$ and (\ref{pulb. of alphaE}), the following diagram commutes
$$\SelectTips{cm}{}\xymatrix{ P\ar[ddd]_{\alpha_P} \ar[dr]^{u}\ar[ddr]_{v}\ar[rr]^{\alpha_P}&&P\ar[ddd]^{u}\\
    &E\ar[ddr]^{\alpha_E}&\!\!\!\!\!\!\!\!\!\!\!\!\!\!\!\!\!\!\!\!\!\!\!\!\shortparallel\\
    \shortparallel\!\!\!\!\!\!\!\!\!\!\!\!\!\!\!\!\!\!\!\!\!\!\!\!&E\ar[dr]_{\alpha_E}&\\
    P\ar[rr]_{v}&&E.}$$
 As $e$ is the  equalizer of $u$ and
$v$ there exists  a morphism $w : P\rightarrow B$ such that
$\alpha_P = ew.$ Therefore, $\alpha_P$ factors through the monic
$e$, hence  $\alpha_\Omega$-dense. $\qquad\square$

In the  topos $M$-{\bf Sets}, the monoid of natural endomorphisms of $\id_{M-{\bf Sets}}$ is isomorphic to the center of $M$. Because corresponding to any natural endomorphism $\alpha: \id_{M-{\bf Sets}} \rightarrow \id_{M-{\bf Sets}}$ the element $\alpha_M(1)=m$ is a central element, i.e., it commutes with all elements of $M$. Conversely, for any central element $m$ of $M$ the mapping $a\mapsto a m$ defines an endomorphism $(\alpha_{m})_A$ of
an arbitrary $M$-set $A$, which is clearly natural with respect to arbitrary action preserving maps.
% (see also, \cite[Example 1.4]{jhonstone1}).

Now in the following we are characterizing all weak topologies of the form $\alpha_{\Omega_M}$ on $M$-{\bf Sets}.
\begin{thm}
    Weak topologies of the form $\alpha_{\Omega_M}$ on $M$-{\bf Sets} where $\alpha$ is a natural endomorphism of $\id_{M-{\bf Sets}}$ are in one to one correspondence to weak ideal topologies $j^{mM}$ on  $M$-{\bf Sets} where $m$ is in the center of $M$.
\end{thm}
{\bf  Proof.} First of all, by the aforesaid paragraph before the Theorem, any  natural endomorphism $\alpha$ of $\id_{M-{\bf Sets}}$ gives a weak ideal topology   $j^{mM}$ on  $M$-{\bf Sets} in which $\alpha_M(1)=m$ is an element in the center of M and conversely, any weak ideal topology of the form $j^{mM}$ for a central element $m\in M$ yields the natural endomorphism $\alpha_{m}$ of $\id_{M-{\bf Sets}}$. Now, we show that these two operations are mutually inverse.
% First, consider a natural endomorphism $\alpha$ of $\id_{M-{\bf Sets}}$. Denote $\alpha_M(1)$ by $m$.
Since the monoid of natural endomorphisms of $\id_{M-{\bf Sets}}$ is isomorphic to the center of M, so it is clear that the assignment $\alpha\mapsto j^{mM}\mapsto \alpha_m$ is identity, i.e. $\alpha_m = \alpha .$ Conversely, recall from~(\ref{def. of jI}) that for any  weak ideal topology  $j^{mM}$ on  $M$-{\bf Sets} in which $m$ is in the center of $M$  we have $(\alpha_m)_\Omega (K) = \{ n\in M|~ nm\in K\} = j^{mM} (K)$ for any right ideal $K$ of $M$. Therefore, $(\alpha_m)_\Omega = j^{mM}$ which shows the assignment $ j^{mM}\mapsto\alpha_m\mapsto $ is identity, i.e., $(\alpha_m)_\Omega =  j^{mM}$.
$\qquad\square$

Recall~\cite{Dikranjan} that a pair $(F, \gamma)$ with an endofunctor
$F: \E \rightarrow \E$ and a natural transformation $\gamma : \id_\E \rightarrow F$ a is called {\it pointed endofunctor} of $\E$. A pointed endofunctor
$(F, \gamma)$ is called a {\it prereflection} if for every commutative diagram
\begin{eqnarray*}
    \SelectTips{cm}{}\xymatrix{ A  \ar[d]_{f} \ar[r]^{\gamma_A}&
FA \ar[d]^{h} \\ B \ar[r]_{\alpha_E } & FB.  }
\end{eqnarray*}
in $\E$ one has $h = F(f)$.

In this regard we have the following remark:
\begin{rk}
{\rm        Let $\alpha$ be a natural endomorphism of the identity
functor on $\E$.   That the pointed endofunctor $(\id_\E, \alpha)$
of $\E$ is a prereflection it follows that the full subcategory
${\rm Fix}(\id_\E, \alpha)$ of $\E$, consisting of all objects
$A\in \E$ such that $\alpha_A$ is isomorphism, is closed under all
(existing) limits of $\E$. In particular, it is replete and closed
under retracts. Furthermore, since $(\id_\E, \alpha)$ is
idempotent, i.e. $\id_\E \alpha = \alpha \id_\E$, so for an object
$A\in \E$, if  $\alpha_A$  be a section, then $A \in {\rm
Fix}(\id_\E, \alpha)$.  (For details see,~\cite[p. 109 and
111]{Dikranjan}.)}
\end{rk}
%%%%%%%%%%%%%%%%%%%%%%%%%%%%%%%%%%%%%%%%%%%%%%%%%%%%%%%%%%%%%%%%%%%%%%%%
\section{The lattice of (productive) weak topologies}
%%%%%%%%%%%%%%%%%%%%%%%%%%%%%%%%%%%%%%%%%%%%%%%%%%%%%%%%%%%%%%%%%%%%%%%%
In this section, we are dealing with three lattices consisting of
weak topologies, productive weak topologies and topologies on a
topos $\E$. We shall describe  their structures in detail.
Finally, in some special toposes, the smallest topology containing a
(productive) weak topology is obtained.

First of all, notice that  the set ${\rm Hom}_{\mathcal{E}}(\Omega , \Omega)$ inherits a partial order from the internal order
on $\Omega$. Indeed, for two morphisms  $j,
k: \Omega \rightarrow \Omega$ of $\E$, one has  $j \leq k$ if and only if $j
= j \wedge k$ where $j \wedge k$ is the composite
$\Omega\stackrel{(j, k)}{\longrightarrow}\Omega\times
\Omega\stackrel{\wedge~}{\longrightarrow}\Omega$.
We  denote by
${\rm Top}(\E),$ ${\rm WTop}(\E)$ and ${\rm PWTop}(\E)$ for the
sets of topologies,  weak topologies and productive weak
topologies on $\E$, respectively. It is clear that these  three sets are subposets of
the poset ${\rm Hom}_{\mathcal{E}}(\Omega , \Omega)$ and we have
\[{\rm  Top}(\E)\subseteq {\rm PWTop}(\E) \subseteq {\rm WTop}(\E). \]
We remark that all these posets have the same binary meets which are
pointwise,  and also they have the top and bottom elements which
are ${\rm true}\circ !_\Omega$ and $\id_\Omega$, respectively.

It is straightforward to verify that the  binary joins   in ${\rm
    WTop}(\E)$ are constructed pointwise, that is, if $j$ and $k$ are
weak topologies, so is the composite $\Omega\stackrel{(j,
    k)}{\longrightarrow}\Omega\times
\Omega\stackrel{\vee~}{\longrightarrow}\Omega$. As $\Omega$ is an
internal distributive lattice in $\E$, it follows that  ${\rm
    WTop}(\E)$ is distributive. If moreover $\E$ is  cocomplete, we can define
\[ (j_1\Rightarrow j_2) = \bigvee\{j\in {\rm WTop}(\E)\mid j\wedge
j_1\leq j_2\},\]
for all  weak topologies $j_1$ and $j_2$ on $\E$,
where $\bigvee$ is the join in the internal complete Heyting
algebra $\Omega$. Note that  the
internal Heyting algebra structure of  $\Omega$ comes from that of
external Heyting algebra structure on the set ${\rm Sub}_\E (E)$
of subobjects of any given object $E$ of $\E$ and  when all colimits exist in $\E$, the complete Heyting algebra structure on the set ${\rm Sub}_\E (E)$ is constructed in the formula (19) in~\cite[p. 497]{maclane}. In this case, the
structure $({\rm WTop}(\E), \vee, \wedge, \Rightarrow, \id_\Omega,
{\rm true}\circ !_\Omega)$ is a complete  Heyting algebra. 

In connection with  the productive weak topologies, one can easily
check that the poset ${\rm PWTop}(\E)$ is a dcpo, i.e. it  has
directed joins which are computed  pointwise. Furthermore, for a
(co)complete topos $\E$, the poset ${\rm PWTop}(\E)$ is a complete
lattice because it has all meets which are calculated pointwise.

On the other hand, as we have already mentioned in Example \ref{weak ideal topology}(i), it is clear that $({\rm
    WTop}(\E), \circ, \id_\Omega)$ is a monoid in which $\circ$ is the
composition of  weak topologies on $\E$, and ${\rm PWTop}(\E)$ is a submonoid of ${\rm WTop}(\E)$.

Let $\E$ be a complete topos. We define two binary
operations $\backslash$ and $/$ on ${\rm WTop}(\E)$  given by
\[
j\backslash k = \bigwedge\{j'\mid j'\in {\rm WTop}(\E), \ j\circ
j'\geq k\},\]
and
\[ k/ j = \bigwedge \{j'\mid j'\in {\rm WTop}(\E), \ j'\circ j\geq
k\}, \]   for  weak topologies $j$ and $k$ on $\E$. We should remark that since meets in ${\rm WTop}(\E)$ are inherited by ${\rm Hom}_{\mathcal{E}}(\Omega , \Omega)$ which is isomorphic to ${\rm Sub}_{\mathcal{E}}(\Omega)$, thus $k/ j$ and $j\backslash k$  also exist whenever $\E$  is a cocomplete topos by the formula (19) in~\cite[p. 497]{maclane}.

Furthermore, it is easily seen that we have
\begin{equation}\label{residuated three}
j\circ j'\geq k \Longleftrightarrow j\geq k/j' \Longleftrightarrow
j'\geq j\backslash k.
\end{equation}
Let us  consider the two posets ${\rm PWTop}(\E)$ and ${\rm WTop}(\E)$ in which the order of each one is opposite, i.e., $\geq$. Then, we have the following (see also, \cite[p. 323, 325]{bir} and
\cite{Jip}):
\begin{thm}\label{lattice of ideals f}
    Let $\E$ be a {\rm(}co{\rm)}complete topos. Then\\
    {\rm (i)} The semilattice $({\rm PWTop}(\E), \wedge)$   is an l-monoid, \\
    {\rm (ii)} The structure $({\rm WTop}(\E), \wedge, \vee, \circ,
    \id_\Omega, \backslash, /)$ is a complete resituated lattice.
\end{thm}
{\bf Proof.}  (i) For productive weak topologies  $j, j'$ and $k$
on $\E$, we have
$$j\circ (j' \wedge k) = j\circ
\wedge \circ (j' , k) = \wedge \circ (j\times j) \circ (j' , k) =
jj'\wedge jk$$
and
$$(j'\wedge k) \circ j
= \wedge \circ (j', k) \circ j = \wedge \circ (j'j, kj) = j'j\wedge kj.$$
This proves (i).\\
(ii)  The assertion is true by (i) and the equivalences mentioned in (\ref{residuated three}). $\quad\quad \square$

For a (co)complete topos $\mathcal{E}$, the inclusion functor \[G_1 :({\rm Top}(\E) , \leq )\rightarrowtail ({\rm WTop}(\E) , \leq) \ ({\rm or \ } G_2 : ({\rm Top}(\E) , \leq )\rightarrowtail ({\rm PWTop}(\E) , \leq)) \]
has a left adjoint, \[F_1 : ({\rm WTop}(\E) , \leq )\rightarrow ({\rm Top}(\E) , \leq ) \ ({\rm or \ } F_2 : ({\rm PWTop}(\E) , \leq )\rightarrow ({\rm Top}(\E) , \leq))\]
which, as any left adjoint to an inclusion, assigns to each
(productive) weak topology $j$ on $\E$ the least topology $j^{\prime}$ on $\E$ with the property
$j\leq j^{\prime}$. We call it the {\it topological reflection} of
$j$ (or {\it idempotent hull} of $j$ in the sense of Dikranjan and Tholen~\cite[p. 82-83]{Dikranjan}). Indeed, we have
\[
j' = \bigwedge \{k\in {\rm Top}(\E)\mid j\leq k\}.
\]
Now we can obtain the following corollary:
\begin{cor}
    The join of a  set of topologies $\{j_\alpha\}_{\alpha\in \Lambda}$ on a {\rm(}co{\rm)}complete topos  $\E$ is the topological reflection of its
    join in ${\rm WTop}(\E)$, i.e. $( \bigvee_{\alpha\in \Lambda}G_1(j_\alpha))^{\prime}$.
\end{cor}
Next we intend to compute the topological reflection of any (productive) weak topology  by a different method.
%Given a (productive) weak topology $j$, one may ask for the
%smallest topology $j'$ containing $j$.  We give a partial
%solution  to this fact here.
Let $\E$ be a complete or cocomplete topos.
Choose a (productive) weak topology $j$ on $\mathcal{E}$. It is
convenient to see that for each natural number $n$, any
$j^n=j\circ\cdots\circ j$ ($n$-times)
is also  a (productive) weak   topology on  $\mathcal{E}$.
Now, one defines the ascending extended ordinal chain of $j$: $$j\leq j^2
\leq j^3\leq \ldots\leq j^\alpha\leq j^{\alpha+1}\leq \ldots\leq
j^\infty\leq j^{\infty+1}$$ as follows: $$j^{\alpha+1} = j\circ
j^\alpha,~~~ j^\beta = \bigvee_{\gamma<\beta}j^\gamma$$ for every
(small) ordinal number $\alpha$ and for $\alpha =\infty$, and for
every limit ordinal $\beta$ and for $\beta =\infty$; here $\infty,
\infty+1$ are (new) elements with $\infty+1>\infty>\alpha$ for all
$\alpha\in {\rm Ord}$, the class of small ordinals.

%The following gives the smallest topology containing  $j$ which is often called {\it
%   idempotent hull} of $j$.
\begin{pro}
Let $j$ be a {\rm(}productive{\rm)} weak topology  on a complete or cocomplete
 topos $\mathcal{E}$. Then $j^\infty$ is the topological reflection of  $j$. That is
 $j^\infty = j'.$
\end{pro}
{\bf Proof.}  It is straightforward
by the corollary of~\cite[p. 83]{Dikranjan}.
$\quad\quad \square$

Now, consider the associated sheaf functor ${\bf a} : \E
\rightarrow {\bf Sh}_{j^\infty}\mathcal{E}$ which is the left adjoint to the inclusion functor  $\iota :  {\bf Sh}_{j^\infty}\mathcal{E}\hookrightarrow \E$. It is  trivial that the restriction functor ${\bf a}$ to the full subcategory  ${\bf
    Sh}_j\mathcal{E}$ of $\E$ and the inclusion functor ${\bf
    Sh}_{j^\infty}\mathcal{E}\hookrightarrow {\bf
    Sh}_j\mathcal{E}$ constitute a geometric morphism between two
toposes ${\bf Sh}_{j^\infty}\mathcal{E}$ and ${\bf
    Sh}_j\mathcal{E}$, for a productive weak topology $j$ on $\E$.

The following result shows that the topos of sheaves associated to the composite of two (weak) topologies can be rewritten as the intersection of two  sheaf toposes.
\begin{thm}
    Let $j$ and $k$ be two weak topologies on a topos $\E$. Then, we have
    $ {\bf Sh}_{jk}(\mathcal{E}) = {\bf Sh}_{j}(\mathcal{E})\cap {\bf Sh}_{k}(\mathcal{E}) = {\bf Sh}_{kj}(\mathcal{E})$ and $ {\bf Sep}_{jk}(\mathcal{E}) = {\bf Sep}_{j}(\mathcal{E})\cap {\bf Sep}_{k}(\mathcal{E}) = {\bf Sep}_{kj}(\mathcal{E}).$ Moreover, if $\mathcal{E}$ is  {\rm(}co{\rm)}complete  then   one has ${\bf Sh}_{j^\infty}\mathcal{E}=
    \bigcap_{\gamma<\infty}{\bf Sh}_{j^\gamma}\mathcal{E}$, as full
    subcategories of $\E$.
\end{thm}
{\bf Proof.}  We only show that ${\bf Sh}_{jk}(\mathcal{E}) = {\bf Sh}_{j}(\mathcal{E})\cap {\bf Sh}_{k}(\mathcal{E}),$ the second assertion is similar.  First of all, we remark that by~\cite[p. 73]{Dikranjan} and the natural isomorphism ${\rm Hom}_{\mathcal{E}}(E, \Omega)\cong {\rm Sub}_{\mathcal{E}}(E)$ (in $E$) for any object $E$ of $\E$, we have $\overline{(\cdot)}^{jk} = \overline{(\overline{(\cdot)}^{k})}^j$ in which $\overline{(\cdot)}^{jk}$ stands for the modal closure operator associated to the composite weak topology $jk$.  Let us now assume that $B\in \E$ be a $jk$-sheaf. We show that $B$ is a $j$-sheaf also. To do so, we first  prove that  any $j$-dense monomorphism $A\rightarrowtail E$ is  $jk$-dense. We have $A\subseteq \overline{A}^{k}$. Then $\overline{A}^{j} \subseteq \overline{(\overline{A}^{k})}^j \subseteq E$. Since $E = \overline{A}^{j}$,  so   $\overline{(\overline{A}^{k})}^j = E$. Thus,  $B$ is a $j$-sheaf. On the other hand, any $k$-dense monomorphism $A\rightarrowtail E$ is  $jk$-dense too. Because if $\overline{A}^{k} = E$ then it is clear that  $\overline{(\overline{A}^{k})}^j = E$. Therefore, $B$ is a $k$-sheaf.

Conversely, we can factor any $jk$-dense monomorphism
$A\rightarrowtail E$  as the composite morphism $A\rightarrowtail
\overline{A}^{k}\rightarrowtail E$ in which $A\rightarrowtail
\overline{A}^{k}$ is a $k$-dense momomorphism and
$\overline{A}^{k}\rightarrowtail E$, a $j$-dense momomorphism. Now
it is evident by the definition of a sheaf that any object $B$ of
$\E$ which is both $j$-sheaf and $k$-sheaf is also a $jk$-sheaf.

%Now we show the last assertion.
Since any $j$-dense monomorphism is clearly
$j^n$-dense,  for any natural number $n$,  we can establish the
chain ${\bf Sh}_j\mathcal{E}\supseteq {\bf    Sh}_{j^2}\mathcal{E}\supseteq \ldots$ as subcategories of $\E$.
Then, ${\bf Sh}_{j^\infty}\mathcal{E}=
\bigcap_{\gamma<\infty}{\bf Sh}_{j^\gamma}\mathcal{E}$ (see also \cite[Corollary A.4.5.16]{jhonstone}).
$\quad\quad \square$
%%%%%%%%%%%%%%%%%%%%%%%%%%%%%%%%%%%%%%%%%%%%%%%%%%%%%%%%%%%%%%%%%%%%%%%%%
\section{Separated  objects}
%%%%%%%%%%%%%%%%%%%%%%%%%%%%%%%%%%%%%%%%%%%%%%%%%%%%%%%%%%%%%%%%%%%%%%%%%
In the present section, we turn our attention to achieve some
necessary and sufficient conditions for the existence of separated object associated to an object of  $\E$ with respect to a  productive weak topology
on $\mathcal{E}$. Afterwards,  among other things, for a
productive weak topology $j$ on  $\E$,  we construct a left
adjoint to the inclusion functor from the category ${\bf
Sep}_j\mathcal{E}$ to the full subcategory of $\E$ consisting of
all objects $E$ of $\E$ for which the closure of diagonal
subobject $\bigtriangleup_E$ of $E\times E$ is closed.

Following~\cite[Corollary V.3.6]{maclane}, for an LT-topology $j$ on a topos $\E$, an object $E$ of  $\E$ and the
diagonal
 $\bigtriangleup_E\in {\rm Sub}_{\mathcal{E}}(E\times E)$,
 the separated object associated to $E$ stands for the coequalizer
\begin{eqnarray}\label{coequalizer}
\SelectTips{cm}{} \xymatrix{\overline{\bigtriangleup_E}\ar@<1.0ex>[r]^{\pi_1}
\ar@<-0.25ex>[r]_{\pi_2} &E\ar@{->>}[r]^{\theta_E} &E',}
\end{eqnarray}
in which $\pi_1$ and $\pi_2$ are the first and second projections.

Note that the following example indicates that for a weak topology  $j$
on a topos $\E$, the object  $E'$ as in (\ref{coequalizer})  is
not separated  in $\E$.
\begin{ex}
{\rm According to  Example \ref{cou. example about sepa.}, it is
straightforward to see that for any $M$-set $E$ the $j^S$-closure
 of $\bigtriangleup_E$ is the following congruence
$$\overline{\bigtriangleup_E} = \{(a , b)\in E\times E\mid \forall s\in S,~ as = bs\}.$$
We  observe that the quotient $M$-set
$M/\overline{\bigtriangleup_M}$ is not separated.
%,i.e. $M' = M/\overline{\bigtriangleup}_M\not \in {\bf Sep}_{j^S}M\!-{\bf Sets}$.
To achieve this
conclusion, fix a non-zero element $s_0\in S$. Also,  consider the
$j^S$-dense monomorphism $SM\subseteq M$ and two action preserving
maps
\[f , g : M\rightarrow M/\overline{\bigtriangleup_M}\] given by
$f(1) = [1]$ and $g(1) = [s_0]$. Since we have $s_0 = 1s_0 \not = s_0s_0 = 0$,
 it yields that $[1] \not = [s_0]$ and hence, $f\not = g$.
 But one can easily check that two restriction maps $f|_{SM}$ and $g|_{SM}$ are equal. We are done.}
\end{ex}
Next, let $j$ be a weak topology  on $\E$ and  $E$ an object of
$\E$.  We assume that $E\stackrel{\theta_E}{\twoheadrightarrow}
S_E\stackrel{\omega_E}{\rightarrowtail}\im (j)^E$ be the image
factorization of the compound morphism $r^E \{\cdot\}_E $ shown in
the diagram below
\begin{equation}\label{epi-mono fact. of r}
\SelectTips{cm}{}\xymatrix{E\ar@{->>}[dr]_{\theta_E} \ar@{->}[rr]^{r^E\{\cdot\}_E}
     &&\im (j)^E\\
&S_E\ar@{>->}[ur]_{w_E} &}
\end{equation}
 in which $r$ is already determined as in
Proposition \ref{im j and idempotent} and $\{\cdot\}_E :
E\rightarrowtail \Omega^E$ stands for the transpose
 of the characteristic map of $\bigtriangleup_E$ denoted as in
 $\delta_E : E\times E \rightarrow \Omega$ (for details, see~\cite{maclane}). Likewise,
 construct the following pullback square in $\E$
 \begin{eqnarray}\label{separated object pullback}
\SelectTips{cm}{}\xymatrix{ E'  \ar@{>->}[d] _{q} \ar@{>->}[r]^{p}&
\Omega_j^E \ar@{>->}[d]^{r^Em^E} \\ S_E \ar@{>->}[r]^{\omega_E ~~} & \im (j)^E.  }
\end{eqnarray}
If $j$ is productive, $E'$ will be separated  as it is a subobject of the $j$-sheaf $\Omega_j^E$.
A reviewing of~\cite{maclane} indicates that for an
LT-topology  $j$, the object $E', $ as in  (\ref{separated object pullback}),
is exactly $S_E$ as well as it  just stands for the coequalizer
as in (\ref{coequalizer}).

The following provides some necessary and sufficient conditions for identifying $E'$ and
$S_E$ as in  (\ref{separated object pullback}).
\begin{thm}\label{A prime and T-A}
 Let $j$ be a weak topology     on $\E$ and  $E$ an object of $\E$. Furthermore,
 let  $w: L_j\rightarrowtail\Omega$   be the equalizer of the morphisms
  $j^2$ and $j$, and $E'$ and $S_E$ as in  the diagram {\rm (}\ref{separated object pullback}{\rm )}. Then the following are equivalent:\\
{\rm (i)}   $E'$ is isomorphic to $S_E$;\\
{\rm (ii)} the subobject $\overline{\bigtriangleup_E}$ of $E\times E$ is closed;\\
{\rm (iii)}  there exists a unique morphism $u: E\times
E\rightarrow L_j$ such that $wu = \delta_E$.
\end{thm}
{\bf  Proof.} (i) $\Longrightarrow$ (ii). Consider the morphism
$k: \im (j)\rightarrowtail \Omega$ as in Proposition \ref{im j and
idempotent} and
 the morphism $m: \Omega_j\rightarrowtail \Omega$ as in the equalizer (\ref{equalizer}) in which
$j= kr$ and $j m= m$. Using (i), by replacing $E'$ by $S_E$  in the pullback diagram (\ref{separated object pullback}), one has
$$\begin{array}{rcl}
w_E\theta_E = r^Em^Ep\theta_E & \Longrightarrow &  r^E \{\cdot\}_E = r^Em^Ep\theta_E\\
& \Longrightarrow & k^Er^E \{\cdot\}_E = k^E r^Em^Ep\theta_E\\
& \Longrightarrow &  j^E \{\cdot\}_E = j^Em^Ep\theta_E\\
& \Longrightarrow &  j^E \{\cdot\}_E = m^Ep\theta_E\\
& \Longrightarrow &  j \delta_E = m{\rm ev}_E(p\theta_E \times \id_E),~~(\textrm{by transposing})
\end{array}$$
in which ${\rm ev}_E$ is the evaluation map on $E$. The last equality
shows that the morphism $j \delta_E$ factors through $m$, the
equalizer (\ref{equalizer}), and then  $j^2 \delta_E = j
\delta_E$. Since ${\rm char} (\overline{\bigtriangleup_E}) = j
\delta_E,  {\rm char}(\overline{\overline{\bigtriangleup_E}})=j^2 \delta_E$ hence, $\overline{\bigtriangleup_E}= \overline{\overline{\bigtriangleup_E}}$,
that is,   $\overline{\bigtriangleup_E}$ is closed in $E\times E$.

(ii) $\Longrightarrow$ (i). By (ii), one has $\overline{\overline{\bigtriangleup_E}}=\overline{\bigtriangleup_E}$ and then
 $j^2 \delta_E = j \delta_E$. In the equalizer (\ref{equalizer}) there exists a unique morphism
 $u : E\times E\rightarrow \Omega_j$ such that
$mu = j \delta_E$. Transposing the last identity gives $m^E \widehat{u} = j^E \{\cdot\}_E$ in which
$\widehat{u} : E\rightarrow \Omega_j^E$ is the transpose of $u$. It follows that
$$\begin{array}{rcl}
 k^E r^Em^E \widehat{u} & = &  j^E m^E\widehat{u}\\
& = &  m^E \widehat{u}\\
& = &  j^E \{\cdot\}_E\\
& = &  k^E r^E \{\cdot\}_E\\
& = &  k^E w_E\theta_E .
\end{array}$$
Since  $k$ is monic,  so is $k^E$ and this yields that
$r^Em^E\widehat{u} = w_E\theta_E$. Then corresponding to
$\widehat{u}$ and $\theta_E$ the pullback situation
(\ref{separated object pullback}) gives a unique morphism $v :
E\rightarrow E'$ such that
 $qv = \theta_E$ and $pv = \widehat{u}$. Since $\theta_E$ is epic the equality $qv = \theta_E$ shows that $q$
 is also  an epic and hence   an isomorphism.

(ii) $\Longleftrightarrow$ (iii). Notice that (ii)  is equivalent to
$j^2 \delta_E = j \delta_E$. Now the desired  equivalence follows from the definition of $L_j$.  $\qquad\square$

We note that a weak topology $j$ on a topos $\E$ need not always
identify $E'$ and $S_E$.
% We note that for a weak  topology $j$ on a topos $\E$ it is not always true   to %identify  $E'$ and $S_E$.
To illustrate   this fact, we give a counterexample through
Theorem \ref{A prime and T-A} as follows.
\begin{ex}\label{s is not sep.}
{\rm According to  Example \ref{cou. example about sepa.},
for any $M$-set $E$ the $j^S$-closure of $\overline{\bigtriangleup_E}$ is the
 following congruence
$$ \overline{\overline{\bigtriangleup_E}} = \{(a , b)\in E\times E\mid \forall s\in S,~ \forall t\in S, ~ast = bst\}.$$ In particular, one has
$\overline{\bigtriangleup_M}\not = \overline{\overline{\bigtriangleup_M}}$, and so $M'\not \cong S_M$. For $ 1S^2 = 0 = s_0S^2$
implies that   $(1, s_0)\in \overline{\overline{\bigtriangleup_M}}$. In contrast,  $s_0 = 1s_0 \not = s_0s_0 = 0$ shows that
  $(1, s_0)\not\in \overline{\bigtriangleup_M}$. However, it is straightforward to see that $\overline{\bigtriangleup_S} =
  \overline{\overline{\bigtriangleup_S}}$ and so, $S'\cong S_S$. Here, $S$ is not $j^S$-separated.
  To achieve this,  fix two distinct
   elements $s ,  t$ in $S$. Moreover, consider the $j^S$-dense monomorphism $SM\subseteq M$
   and two  action preserving  maps $f , g : M\rightarrow S$
   defined by $f(1) = s$ and $g(1) = t$. Then one has $f\not = g$, but $f|_{SM} = 0 = g|_{SM}$.}
\end{ex}
Recall that the kernel pair of a map $t : E\rightarrow W$
is the pullback of the diagonal $\bigtriangleup_W$ along $t\times t$,
 as in the following pullback situation
 \begin{eqnarray*}
\SelectTips{cm}{}\xymatrix{B \ar@{>->}[d]_{} \ar[r]& \bigtriangleup_W
\ar@{>->}[d] \\ E\times E \ar[r]^{t\times t} & W\times W.}
\end{eqnarray*}
The following shows that if  we identify $E'$ and $S_E$, the
morphism $\theta_E$ is as close to being a monic.
\begin{lem}\label{theta-A and its kernel pair}
Let $j$ be a productive weak topology     on $\E$ and  $E$ an object of $\E.$ If we identify $E'$ and $S_E$ defined as in
 {\rm (}\ref{separated object pullback}{\rm )}, then the kernel pair of $\theta_E : E\twoheadrightarrow S_E$ is precisely the closure
  $\overline{\bigtriangleup_E}$ of the diagonal $\bigtriangleup_E\subseteq E\times E$.
\end{lem}
{\bf  Proof.} Since $E'$ and $S_E$ are the same, the following
diagram stands for the image factorization of the compound
morphism $r^E \{\cdot\}_E$,
 \begin{eqnarray}\label{image of r-A}
\xymatrix{  E\ar@{->>}[d] _{\theta_E} \ar[r]^{\{\cdot\}_E}& \Omega^E \ar[d]^{r^E} \\ S_E ~\ar@{>->}[r]_{r^Em^Ep~~~~} & \im (j)^E}
\end{eqnarray}
in which $w_E = r^Em^Ep$ in the diagram~(\ref{separated object pullback}). It yields that
the kernel pair of $\theta_E$ is the same as that of $r^E \{\cdot\}_E$. Suppose now that $(f, g)$
is the kernel pair of  $r^E \{\cdot\}_E$. One has $r^E \{\cdot\}_E f = r^E
\{\cdot\}_E g$.
Transposing this equality first and then  composing with $k$  gives
$j \delta_E (f\times \id_E) = j \delta_E (g\times \id_E)$, since
$j = kr$. Similar to the proof of~\cite[Lemma V.3.5]{maclane},
we observe that $(f, g)$  is contained in  $\overline{\bigtriangleup_E}$.
 By Lemma~\ref{being j- separated of omega-j} the object $\Omega_j$ is a sheaf.
 Also, we can deduce that
  $r^Em^Ep\theta_E\pi_1 = r^Em^Ep\theta_E\pi_2$ and then by (\ref{image of r-A}),
  $r^E\{\cdot\}_E\pi_1 = r^E\{\cdot\}_E\pi_2$. This implies that
   $\overline{\bigtriangleup_E}$  is contained in $(f, g)$.   The proof is now complete. $\qquad\square$

Now we will focus on establishing some transfer property from LT-topologies  on a topos $\E,$ which can be found in
\cite[Lemma V.3.3]{maclane},
 to weak topologies on $\E$.
\begin{lem}\label{delta is closed}
Let  $j$ be a weak   topology on a topos $\mathcal{E}$. Then for any object $E$ of  $\mathcal{E}$ the following are equivalent:\\
{\rm (i)}  $E$ is separated;\\
{\rm (ii)} the diagonal $\bigtriangleup_E\in {\rm Sub}_{\mathcal{E}}(E\times E)$ is a closed subobject of $E\times E$;\\
{\rm (iii)}  $j^E\circ \{\cdot\}_E = \{\cdot\}_E$, as in the commutative diagram
$$\xymatrix{E \ar[dr]_{\{ \cdot\}_E} \ar[r]^{\{ \cdot\}_E} & \Omega ^E \ar[d]^{j^E} &\\ &\Omega ^E; &}$$ {\rm (iv)} for any $f : A \rightarrow E$, the graph of $f$ is a closed subobject of $A\times E$.
\end{lem}
{\bf  Proof.} The parts (i) $\Longrightarrow$ (ii) $\Longrightarrow$ (iii) $\Longrightarrow$ (iv) are proved in a similar way to
\cite{maclane}.  Also, the proof of (ii) $\Longrightarrow$ (i)  is analogous to  Proposition 9.2.4
 of~\cite[Vol. III]{handbook1}. Finally, (iv) $\Longrightarrow$ (ii) is established by setting $f = \id_E$, as  $\bigtriangleup_E$
 stands for the graph of $\id_E$.   $\qquad\square$

For a  productive weak topology $j$ on a topos $\E$, we write
$\mathcal{C}_j$ for the full subcategory of $\E$ consisting of all
objects $E$ of $\E$ for which the subobject
$\overline{\bigtriangleup_E}$ of $E\times E$ is closed, i.e. $\overline{\bigtriangleup_E}=\overline{\overline{\bigtriangleup_E}}$.  It is
immediate that $1\in \mathcal{C}_j$ and $\mathcal{C}_j$ is closed
under finite products.  Notice that in Example~\ref{s is not sep.}
the semigroup $S$ (as a sub $M$-set of $M$) is not separated,
however, $\overline{\bigtriangleup_S}$ is closed in $S\times S$
with respect to the weak topology $j^S$ on $M$-{\bf Sets}. Therefore in general ${\bf
Sep}_j\mathcal{E}$ and $\mathcal{C}_j$ do not necessarily
coincide as subcategories of $\E$.

 In the next remark we provide a short characterization of the former subcategory of $\E$.
\begin{rk} {\rm It is immediate
that for a principal productive weak topology $j$ on a topos $\E$,
one has $E\in \mathcal{C}_j$ if and only if
$\overline{\overline{\bigtriangleup_E}}^{^\circ}\subseteq
\bigtriangleup_E$ as subobjects of $E\times E$.
 Furthermore,, for
a natural endomorphism $\alpha$ of the identity functor on $\E$,
an easy computation shows that a subobject $m : A\rightarrowtail
E$ in $\E$ is $\alpha_\Omega$-closed if and only if one has ${\rm
char}(m) = {\rm char}(m) \alpha_E$. Then we have $A\in
\mathcal{C}_{\alpha_\Omega}$ if and only if $\alpha_\Omega
\delta_A = \delta_A$ if and only if $\delta_A = \delta_A \circ
(\alpha_A \times \alpha_A)$. This means that $A\in
\mathcal{C}_{\alpha_\Omega}$ if and only if the morphism $\alpha_A
: A\rightarrow A$ be monic. Note that for any $E\in \E$, one has
$\overline{\bigtriangleup_E}^{^\circ} =
\overline{\bigtriangleup_E}\circ\im (\alpha_{\overline{E}})$.
Also,  $\overline{\bigtriangleup_E} \circ\im
(\alpha_{\overline{E}}) \subseteq \im (\alpha_{E}) \times \im
(\alpha_{E})$ as subobjects of $E\times E.$}
\end{rk}
Now we construct an adjunction.
\begin{cor}\label{a refle. subcate. c-j}
    For any productive weak topology  $j$ on a topos $\mathcal{E}$,
    the inclusion functor ${\bf Sep}_j\mathcal{E}\rightarrowtail
    \mathcal{C}_j$ has a left adjoint $ L :
    \mathcal{C}_j\longrightarrow {\bf Sep}_j\mathcal{E}$ defined by
    $E\mapsto E'$.
\end{cor}
{\bf Proof.} First of all, by Lemma~\ref{delta is closed} ((i)
$\Leftrightarrow$ (ii))  the inclusion functor  ${\bf
Sep}_j\mathcal{E}\rightarrowtail \mathcal{C}_j$ exists. Let $E$ be
an object of $\mathcal{C}_j$. By Theorem~\ref{A prime and T-A},
$E'$ is isomorphic to $S_E$ and then, Lemma~\ref{theta-A and its
kernel pair} shows that the kernel pair of $\theta_E :
E\twoheadrightarrow E'$ is precisely the closure
$\overline{\bigtriangleup_E}$ of the diagonal
$\bigtriangleup_E\subseteq E\times E$. Now let $F$ be a separated
object and $f: E\arrow F$ a morphism in $\mathcal{C}_j$.  It is
clear that the diagram
$$\SelectTips{cm}{} \xymatrix{\bigtriangleup_E\ar@<1.0ex>[r]^{\pi'_1}
    \ar@<-0.25ex>[r]_{\pi'_2} &E\ar[r]^{f} &F,}$$
commutes in which $\pi'_1$ and $\pi'_2$  are the first and second projections. Then   the following diagram  is also
$$\SelectTips{cm}{} \xymatrix{\overline{\bigtriangleup_E}\ar@<1.0ex>[r]^{\pi_1}
    \ar@<-0.25ex>[r]_{\pi_2} &E\ar[r]^{f} &F}$$
commutative since  $F$ is separated.  From the coequalizer
diagram~(\ref{coequalizer}), we can deduce that there is a unique
morphism $g : E'\arrow F $ such that $g\theta_E =
f$.$\quad\quad\square$

In particular we have the following
\begin{rk}\label{Fryed}
    {\rm    Recall~\cite[p. 87]{Freyd} that if    $\mathcal{B}$ is  a complete, cocomplete and co-well-powered category and $\mathcal{A}$ a full subcategory replete in $\B$ such that $\mathcal{A}$ is closed under the formation of products and subobjects then, $\mathcal{A}$ is a reflective subcategory of $\B$.

        Next let us suppose that $\E$ be  a complete, cocomplete and co-well-powered topos and   $j$ a productive weak  topology on $\mathcal{E}$.
        Since the full subcategory ${\bf Sep}_j\mathcal{E}$ of $\E$ is closed under the formation of products and subobjects, it follows that it is a reflective subcategory of $\mathcal{E}$.
        Thus, the inclusion functor ${\bf Sep}_j\mathcal{E}\rightarrowtail
        \mathcal{E}$ has a left adjoint $ R : \mathcal{E}\longrightarrow {\bf Sep}_j\mathcal{E}$.  One can construct the functor
        $R$ by the adjoint functor theorem~\cite[V. 18.12]{Adamek}.
        For more discussions see   also~\cite[Definition 3.1]{Clementino}.}
\end{rk}
The notion of a {\it quasitopos} can be found in~\cite{jhonstone}.
With the assumption given by the last paragraph we can show that ${\bf Sep}_j\mathcal{E}$ is a quasitopos.  Indeed
% The following indicates the basic fact we are interested in.
\begin{thm}\label{sep j is a quasitopos}
Let $\mathcal{E}$ be a complete, cocomplete and co-well-powered
topos and $j$ a productive weak  topology on $\mathcal{E}$. Then
the category ${\bf Sep}_j\mathcal{E}$ is a quasitopos.
\end{thm}
{\bf  Proof.} As we have already mentioned after Example~\ref{weak ideal topology}, ${\bf
Sep}_j\mathcal{E}$ is a (finitely) complete  category as well as
it is closed under exponentiation. Since ${\bf
Sep}_j\mathcal{E}$ is a reflective subcategory of $\mathcal{E}$
it indicates that, by~\cite[Vol. I, Proposition 3.5.4]{handbook1},
${\bf Sep}_j\mathcal{E}$ is a  cocomplete category. Afterwards,
let $B$ be
 an object of ${\bf Sep}_j\mathcal{E}$ and $j_B =  j\times \id_B,$  the induced weak
 topology by $j$ on the slice topos $\mathcal{E}/B$ (see also~\cite{KMItalian}).
  One can easily check that
\begin{equation*}
{\bf Sep}_{j_B}\mathcal{E}/B \cong {\bf Sep}_j\mathcal{E}/B.
\end{equation*}
From which it follows that the category ${\bf Sep}_{j_B}\mathcal{E}/B$ is cartesian closed. Moreover, any $j$-dense monomorphism
 in ${\bf Sep}_j\mathcal{E}$ is epic, by the definition of a separated object in $\E$. Meanwhile, any strong monomorphism in
  ${\bf Sep}_j\mathcal{E}$ must be $j$-closed. For establishing this, associated to any strong monomorphism $i:C\rightarrow
  E$ in ${\bf Sep}_j\mathcal{E}$ there is a factorization  as follows
$$\SelectTips{cm}{}\xymatrix{  C \ar@{>->}[d] _{\iota} \ar@{>->}[r]^{i}&E  \\ \overline{C}\ar@{>->}[ru]_{\overline{i}}  &}$$
The morphism $\iota$ is $j$-dense and thus it is epic in ${\bf
Sep}_j\mathcal{E}$. Since $i$ is strong it follows that there
exists  a unique morphism $w$ in the  commutative square below
$$\SelectTips{cm}{}\xymatrix{C \ar[r] ^{\id_C} \ar@{>->}[d]_{\iota}&
C \ar@{>->}[d]^{i} \\ \overline{C} \ar[r]_{\overline{i}}\ar@{-->}[ur]|w & E}$$
such that $w\iota = \id_C, iw=\bar{i}$. Now we have
 $\overline{C}\subseteq C.$
This yields that $i$ is $j$-closed. Finally Lemma \ref{being j- separated of omega-j} and the bijection
 (\ref{omega j classifies closed subobjects}) show that $\Omega_j$ is a
 weak subobject classifier for ${\bf Sep}_j\mathcal{E}$.
 $\qquad\square$
%%%%%%%%%%%%%%%%%%%%%%%%%%%%%%%%%%%%%%%%%%%%%%%%%%%%%%%%%%%%%%%%%%%%%%%%%%%%%
\section{An adjunction}
%%%%%%%%%%%%%%%%%%%%%%%%%%%%%%%%%%%%%%%%%%%%%%%%%%%%%%%%%%%%%%%%%%%%%%%%%%%%%
In the previous section we were concerned with separated object associated to any  object of a  topos $\mathcal{E}$.
 Now we shall obtain, among other things,  the sheaf associated to any separated object of $\mathcal{E}$ with respect
 to a given  productive weak  topology $j$ on $\mathcal{E}$.

 Following Lemma \ref{delta is closed}(ii), for any separated object $E$ of  $\E$ the diagonal $\bigtriangleup_E$ is a
  closed subobject of $E\times E$. In this case, the characteristic map of
  $\bigtriangleup_E$ denoted by $\delta_E : E\times E \rightarrow \Omega$ satisfies
 $j\delta_E = \delta_E$. Then    the equalizer diagram (\ref{equalizer})  gives a unique morphism
\begin{equation}\label{arrow t of a sep. obj. A}
\gamma_E : E\times E\rightarrow \Omega_j \ \mbox{such that} \ m  \gamma_E =\delta_E .
\end{equation}
  We denote the exponential transpose of $\gamma_E $ by $\widehat{\gamma}_E  : E\rightarrow \Omega_j^E$; that is
\begin{equation}\label{Transpose of alpha}
ev_E (\widehat{\gamma}_E \times \id_E ) = \gamma_E .
\end{equation}
Indeed in the next lemma we will show that   $\widehat{\gamma}_E $ is monic.
\begin{lem}\label{being mono t-A's}
Let  $j$ be a weak topology     on a topos $\E$ and $E$ a
separated object  of  $\E$. Then the morphism
$\widehat{\gamma}_E $, defined by the transpose of $\gamma_E $ as in
(\ref{arrow t of a sep. obj. A}), is a monomorphism.
\end{lem}
{\bf  Proof.} Consider an object $A$ and two morphisms $g, h :
A\rightarrow E$ for which $\widehat{\gamma}_E  g =
\widehat{\gamma}_E  h$. By taking the exponential transpose of this
equality we obtain $\gamma_E (g\times \id_E) = \gamma_E (h\times
\id_E)$. Now we have
$$\delta_E (g\times \id_E) =  m \gamma_E (g\times \id_E) =m \gamma_E (h\times \id_E)= \delta_E (h\times \id_E) .$$
Using the proof of Lemma IV.1.1 as in~\cite{maclane}, we get $g = h$. $\qquad\square$
\begin{rk}
\rm{As a consequence of the above lemma and Lemma~\ref{being j- separated of omega-j}, we can mention that the separated objects of $\E$ are precisely subobjects of sheaves on $\E$ for a productive weak topology $j$ on $\E$.  Thus two subcategories ${\bf Sep}_j\mathcal{E}$ and ${\bf Sh}_j\mathcal{E}$ of $\E$ coincide if and only if ${\bf Sh}_j\mathcal{E}$ is closed
    under subobjects in $\E$. Furthermore, we can see that every object of $\E$ is
    injective for the class of $j$-dense monomorphisms if and only if every
    $j$-dense monomorphism is split. (For a more general case, see~\cite[Theorem 2.1]{jhonstone1}.)}
\end{rk}

It is well known~\cite{maclane} that for an LT-topology  $j$ on
$\E$,  $\Omega_j = \im (j)$, and for a separated object $E$ of
$\E$,  the {\it sheaf associated} to $E$ is  the closure of the
composite morphism $r^E \{\cdot\}_E : E\rightarrowtail \Omega_j^E$
in which $r$
 is already determined as in Proposition \ref{im j and idempotent}.
Next we provide the relationship between the sheaf associated to
$E$ and the morphism $\gamma_E $ as in (\ref{arrow t of a sep. obj.
A}) in $\E$.
\begin{pro}\label{asso. sheaf on A and t-A}
For an LT-topology  $j$ on $\E$ and a separated object $E$ of  $\E$, we have $\widehat{\gamma}_E  = r^E \{\cdot\}_E$.
\end{pro}
{\bf  Proof.} To check our claim, we write
 \begin{eqnarray}\label{t hat and the singelton}
m ev_E (\widehat{\gamma}_E \times \id_E) & =& m\gamma_E \nonumber\\
& =& \delta_E = j  \delta_E\nonumber\\
&=&  j ev_E (\{\cdot\}_E\times \id_E).
\end{eqnarray}
The exponential transpose of (\ref{t hat and the singelton}) yields that $m^E\widehat{\gamma}_E  = j^E \{\cdot\}_E$.
But $j^E = m^E r^E$. Then, $m^E\widehat{\gamma}_E  = m^E r^E\{\cdot\}_E$. Since the functor $(-)^E$
preserves monomorphisms,
 $m^E$ is a monomorphism and hence, $\widehat{\gamma}_E  = r^E \{\cdot\}_E$. $\qquad\square$

We remark that in an analogous way to~\cite{Barr,carboni}, we can define a
{\it weak    topology {\rm (}modal closure operator}) on a category
 with finite limits. Hence, for a weak  topology $j$ on a topos $\mathcal{E}$,
 the notion of $j$-sheaves can be defined in the finite
 complete category ${\bf Sep}_j\mathcal{E}$.

The following  determines $j$-sheaves in  ${\bf Sep}_j\mathcal{E}$.
\begin{pro}\label{sheaves in the cat. sep. objects}
Let  $j$ be a productive weak    topology on $\mathcal{E}$ and  $E$ a separated object of  $\E$.
Then the following conditions are equivalent: \\
{\rm (i)}  $E$ is a $j$-sheaf in  ${\bf Sep}_j\mathcal{E}$;\\
{\rm (ii)} $E$ is a $j$-sheaf in  $\mathcal{E}$;\\
{\rm (iii)}  $E$ is closed in $\Omega_j^E$.
\end{pro}
{\bf  Proof.} (i) $\Longrightarrow$ (iii). Corresponding to the
separated object $E$ consider the monomorphism  $\widehat{\gamma}_E  :
E\rightarrowtail \Omega_j^E$  defined as in~(\ref{Transpose of
alpha}).
 The inclusion morphism  $\iota : E\rightarrowtail \overline{E} $ is $j$-dense
 in $\E$ and so is epic in ${\bf Sep}_j\mathcal{E}$. We point out here that
 $\overline{E}$ is separated as it is a subobject of the $j$-sheaf $\Omega_j^E$.
  Since  $E$ is a $j$-sheaf in  ${\bf Sep}_j\mathcal{E}$, we obtain
  a retract $q$ of $\iota$ in the following situation
$$\SelectTips{cm}{}\xymatrix{E \ar@{>->}[d]_{\iota}\ar[r]^{\id_E} & E  \\ \ar@{-->}[ur]_{q}\overline{E} .}$$
In this way, we have $\iota q \iota = \iota\id_E = \id_{\overline{E}}\iota$.
That $\iota$ is epic  yields that $\iota q =  \id_{\overline{E}}$ and so,
 $q$  is an isomorphism. Then we have $E = \overline{E}$, as subobjects of
 $\Omega_j^E$. Therefore, $E$ is closed in $j$-sheaf $\Omega_j^E$.\\
(ii) $\Longrightarrow$ (i)  is clear. \\
(iii) $\Longleftrightarrow$ (ii). An investigation of~\cite[Lemma V.3.4]{maclane} for productive weak topologies on $\mathcal{E}$ instead of LT-topologies,
  shows that a separated object $E$ is closed in $\Omega_j^E$ if and only if $E$ is a $j$-sheaf in  $\mathcal{E}$. $\qquad\square$

Finally, in what follows we provide the sheaf associated to a separated object in a   topos $\mathcal{E}$.
\begin{thm}\label{closure of A is a sheaf}
Let $j$ be  a productive weak    topology   on a topos
$\mathcal{E}$ and  $E$ a separated object  of  $\E$. Moreover, let
$\overline{E}$ be the closure of $E$ as a subobject of
$\Omega_j^E$ via the morphism $\widehat{\gamma}_E $. Then
$\overline{E}$ is a $j$-sheaf in $\mathcal{E}$.
\end{thm}
{\bf  Proof.} Using Proposition \ref{sheaves in the cat. sep. objects},
it is enough to show that $\overline{E}$ is closed in $\Omega_j^E$, i.e.
 $\overline{\overline{E}} = \overline{E}$ as subobjects of  $\Omega_j^E$.
 We denote the closure of $\widehat{\gamma}_E  : E\rightarrowtail \Omega_j^E$ by
 $u : \overline{E}\rightarrowtail \Omega_j^E$ and the closure of $u$ by
 $\overline{u} : \overline{\overline{E}}\rightarrowtail \Omega_j^E$.
  We note that $ \overline{E}$ and $\overline{\overline{E}}$ are separated as they
 are subobjects of the $j$-sheaf $\Omega_j^E$. The morphism
 $\nu : \overline{E}\rightarrowtail \overline{\overline{E}} $ in the
  following commutative  diagram is $j$-dense in $\E$,
\begin{eqnarray}\label{diagram for topologies}
\SelectTips{cm}{}\xymatrix{\overline{E}
\ar@{>->}[d]_{\nu}\ar@{>->}[r]^{u~} & \Omega_j^E \\
  \ar@{>->}[ur]_{\overline{u}}\overline{\overline{E}} }
\end{eqnarray}
One has, $j{\rm char}(\widehat{\gamma}_E ) u = {\rm true}_{\overline{E}}$ as in the following diagram
\begin{eqnarray}\label{diagram for u and t hat}
\SelectTips{cm}{}\xymatrix{\overline{E}\ar@{>->}[dd]_{u}\ar[rrrr]&&&&1\ar[dd]^{{\rm true}}\\
    &E\ar[rr] \ar@{>->}[d]_{\widehat{\gamma}_E } && 1\ar[d]^{{\rm true}}&\\
    \Omega_j^E \ar@{=}[r]&\Omega_j^E  \ar[rr]_{{\rm char}(\widehat{\gamma}_E )}&& \Omega\ar[r]_{j}&\Omega}
\end{eqnarray}
 Then using (\ref{diagram for topologies}) we observe that
 $j{\rm char}(\widehat{\gamma}_E ) \overline{u} \nu = ({\rm true}_{\overline{\overline{E}} })
 \nu$. Since $\nu$ is dense and $\Omega_j^E$ separated, it follows that
 $j{\rm char}(\widehat{\gamma}_E ) \overline{u}= {\rm true}_{\overline{\overline{E}} }$.
 Hence, in the diagram (\ref{diagram for u and t hat})
 there exists   a unique morphism $w : \overline{\overline{E}}\rightarrow \overline{E} $
 such that $u w = \overline{u}$. It follows that $\overline{\overline{E}} = \overline{E}$ as required. $\qquad\square$

%According to Theorem \ref{closure of A is a sheaf}, in a same way
%as in~\cite[Corollary V.3.7]{maclane}, one has:
Now we can deduce easily the following corollary:
\begin{cor}\label{sheaf asso.sep}
For any productive weak topology $j$ on a topos $\mathcal{E}$, the
inclusion functor ${\bf Sh}_j\mathcal{E}\rightarrowtail
 {\bf Sep}_j\mathcal{E}$ has a left adjoint
 $S : {\bf Sep}_j\mathcal{E}\longrightarrow {\bf Sh}_j\mathcal{E}$ defined by $E\mapsto
 \overline{E}$ as subobjects of  $\Omega_j^E$.
\end{cor}
At this point let us sum up our result. For any object $E$ of
$\mathcal{C}_j$ first construct the epimorphism $\theta_E : E\twoheadrightarrow
S_E$ defined as in~(\ref{epi-mono fact. of r}). Then $S_E\cong
E'$. We can embed the separated object $E'$ into $\Omega_j^{E'}$
by the monomorphism  $\widehat{\alpha}_{E'}$ and then we denote
the resulting morphism
\[ E\stackrel{\theta_E}{\twoheadrightarrow} S_E\cong E'\hookrightarrow  \overline{E'}\rightarrowtail \Omega_j^{E'}\] by $i_E$. As $E'\hookrightarrow \overline{E'}$ is a $j$-dense monomorphism in  ${\bf Sep}_j\mathcal{E}$ so it is epi and hence  the composite morphism
$ E\stackrel{\theta_E}{\twoheadrightarrow} S_E\cong E'\hookrightarrow  \overline{E'}$ is an epi too which is precisely the $\im (i_E)$. Then, for any object $E$ of $\mathcal{C}_j$,   we may find  a map $i_E : E\rightarrow I_E$ where $I_E$ is a sheaf and it has  minimal kernel pair $\overline{\bigtriangleup}_E$. Therefore  the compound functor
$ SL :\mathcal{C}_j\longrightarrow {\bf Sh}_j\mathcal{E}$
assigns to any $E$ of $\mathcal{C}_j$ the sheaf $\overline{\im(i_E)}$.

We now  provide the associated sheaf  functor with respect to a  weak topology $j$ on $\mathcal{E}$.
\begin{thm}
    With the notations given in corollary~(\ref{sheaf asso.sep}) the composite  functor $ SL :\mathcal{C}_j\longrightarrow {\bf Sh}_j(\mathcal{E})$ is a left adjoint to  the inclusion functor ${\bf Sh}_j(\mathcal{E})\rightarrowtail \mathcal{C}_j$
    which assigns to any object $E$ of $\mathcal{C}_j$ the sheaf $\overline{\im(i_E)}$ as a subobject of  $I_E$. This functor preserves the terminal object, monics and finite products. In the case of a  complete, cocomplete and co-well-powered
    topos $\mathcal{E}$ the inclusion functor ${\bf Sh}_j(\mathcal{E})\rightarrowtail \mathcal{E}$ has the composite  left adjoint
    $SR : \mathcal{E}\longrightarrow {\bf Sh}_j(\mathcal{E})$.
\end{thm}
{\bf Proof.} By applying Corollaries~\ref{a refle. subcate. c-j} and \ref{sheaf asso.sep} and Remark~\ref{Fryed}, we get the result. $\quad\square$
%%%%%%%%%%%%%%%%%%%%%%%%%%%%%%%%%%%%%%%%%%%%%%%%%%%%%%%%%%%%%%%%%%%%%

\end{document}